\documentclass[english,envcountsect, envcountsame]{svjour3}
\usepackage[T1]{fontenc}
\usepackage[latin9]{inputenc}
\usepackage{array}
\usepackage{bbding}
\usepackage{mathrsfs}
\usepackage{url}
\usepackage{enumitem}
\usepackage{amsmath}
\usepackage{amssymb}
\usepackage{graphicx}
\PassOptionsToPackage{normalem}{ulem}
\usepackage{ulem}

\makeatletter

\providecommand{\tabularnewline}{\\}

\newenvironment{svmultproof}{\begin{proof}}{\qed\end{proof}}

\DeclareMathOperator{\Trim}{Trim}
\DeclareMathOperator{\red}{red}
\DeclareMathOperator{\Co}{Co}
\setlist[enumerate,1]{label = (\alph*)}

\makeatother

\usepackage{babel}

\begin{document}
\title{Invariants for metrisable locally compact Boolean spaces}
\author{Andrew B. Apps\thanks{I gratefully acknowledge the use of the library facilities of the
University of Cambridge.\protect\phantom{\Envelope{}}}}
\institute{\Envelope{}  Andrew Apps, Independent researcher, St Albans, UK\\
$^{\phantom{\Envelope{}}}$ andrew.apps@apps27.co.uk\\\\}
\date{December 2024}
\maketitle
\begin{abstract}
Pierce identified 3 invariants of a compact metrisable Boolean space,
derived from its \emph{Cantor-Bendixson sequence}, that determine
the space up to homeomorphism. For locally compact spaces we define
an additional invariant, the \emph{compact rank}, and show that these
4 invariants determine a locally compact metrisable Boolean space
up to homeomorphism. We also identify which combinations of the 4
invariants can arise in practice.

A Boolean ring and its associated Boolean space are \emph{primitive}
if the ring is disjointly generated by its \emph{pseudo-indecomposable}
(PI) elements. Spaces in this important sub-class of Boolean spaces
can be well described (uniquely in the case of compact spaces) by
an \emph{extended PO system} (poset with a distinguished subset).
We define the Cantor-Bendixson sequence and associated invariants
for a PO system, and show that almost all of the invariant information
for a primitive space can be recovered from that of an associated
extended PO system.

We also show how the primitivity\emph{ }of a Boolean space corresponds
to a notion of primitivity of the \emph{additive measure }associated
with the \emph{rank function }of a space, which in turn depends on
the additive measure being sufficiently ``self-similar''. We use
these ideas to develop a method for constructing non-primitive spaces.
\end{abstract}

\section{Introduction}

Stone's celebrated paper~\cite{Stone} established the duality between
Boolean rings (which may or may not have a multiplicative identity)
and \emph{Boolean spaces}, defined as totally disconnected locally
compact Hausdorff topological spaces. The structure of compact metrisable
Boolean spaces has been studied by authors such as Pierce and Ketonen
using a series of invariants. We extend these results to locally compact
metrisable Boolean spaces, which we term \emph{$\omega$-Stone spaces},
being the Stone spaces of countable Boolean rings, and examine the
interplay of these invariants with other invariants that arise for
\emph{primitive }Boolean spaces.

The isolated point structure of a Boolean space $W$ is captured by
the \emph{Cantor-Bendixson sequence}: let $W^{(1)}=W^{\prime}$ comprise
all the non-isolated points of $W$, and proceed inductively to define
$W^{(\xi)}$ for any ordinal $\xi$. Let $W^{[\xi]}=W^{(\xi+1)}-W^{(\xi)}$;
for increasing~$\xi$, the $\overline{W^{[\xi]}}$ form a decreasing
series of closed subsets of $W$. Three invariants can be defined
based on this sequence, which Pierce showed determine the structure
of a compact metrisable Boolean space up to homeomorphism:
\begin{enumerate}
\item $\nu(W)$, the smallest $\xi$ such that $W^{[\xi]}=\emptyset$;
\item $n(W)$, the size of $W^{[\mu]}$ if $\nu(W)$ is the successor ordinal
$\mu+1$;
\item a \emph{rank function }$r_{W}\colon K(W)\rightarrow\omega_{1}$, the
set of all countable ordinals, where $K(W)=W^{(\nu(W))}$, the \emph{perfect
kernel }of $W$, being its largest subspace $C$ such that $C^{\prime}=C$.
\end{enumerate}
As $K(W)$ has no isolated points, it is either empty or homeomorphic
for metrisable~$W$ to the Cantor set or to the Cantor set minus
a point. A space $W$ is \emph{scattered }if $K(W)=\emptyset$ and
\emph{uniform} if $\overline{W^{[\xi]}}\cap K(W)$ is non-empty for
all $\xi<\nu(W)$. Any compact metrisable Boolean space can be expressed
as the disjoint union of two clopen subsets, one uniform or empty
and the other scattered or empty.

We introduce an additional invariant to capture the compactness structure
of a locally compact $\omega$-Stone space $W$, namely the \emph{compact
rank} $\rho(W)$, which is the smallest $\xi$ such that $\overline{W^{[\xi]}}$
is compact. 

It is well known that a compact scattered $\omega$-Stone space is
homeomorphic to the ordinal space $\omega^{\mu}.n+1$, equipped with
the order topology, for some ordinal $\mu$ and $n\in\mathbb{N}$.
Less well known is the result of Flum and Martinez that \uline{any}
scattered $\omega$-Stone space is homeomorphic to some ordinal space.
Specifically, any such space $W$ is homeomorphic to a space of the
form $\omega^{\mu}.n+1$ or $\omega^{\rho}$, or to a disjoint union
of the two with $\rho\leqslant\mu$, where $\mu$ and $\rho$ are
ordinals; $W$ is uniquely determined by the tuple $[\nu(W),\rho(W),n(W)]$. 

In the first part of this paper we extend these existing results by
showing that any $\omega$-Stone space is determined up to homeomorphism
by the tuple

\[
[(K(W),r_{W}),\nu(W),\rho(W),n(W)],
\]

and conversely that for any suitable tuple $[(K,r),\nu,\rho,n]$ there
is an $\omega$-Stone space yielding these invariants. Moreover, any
$\omega$-Stone space can be written as the disjoint union of two
clopen subsets, one strongly uniform or empty and the other scattered
or empty, and strongly uniform spaces are determined up to homeomorphism
by their rank functions. Here, an $\omega$-Stone space is \emph{strongly
uniform} if it can be expressed as a disjoint union of compact open
uniform subspaces.

In the second part of the paper we consider primitive Boolean spaces
and their invariants. Primitive Boolean algebras were introduced by
Hanf~\cite{Hanf} as a class of well-behaved Boolean algebras, and
the definition can be readily extended to Boolean rings which may
or may not have a $1$: we recall that an element $A$ of a Boolean
ring $R$ is \emph{pseudo-indecomposable (PI)} if for all $B\in(A)$,
either $(B)\cong(A)$ or $(A-B)\cong(A)$, where $(A)=\{C\in R\mid C\subseteq A\}$;
and $R$ and its Stone space are \emph{primitive }if every element
of $R$ is the disjoint union of finitely many PI elements. We follow
most subsequent authors after Hanf in dropping the requirement that
$1_{R}$ be PI for a primitive Boolean algebra $R$. For any primitive
Boolean space $W$ there is an associated \emph{extended PO system}
$(P,L,f)$, where $P$ is the canonical \emph{PO system} (poset with
a distinguished subset) of $W$ consisting of the homeomorphism classes
of compact open PI subsets, and a well-behaved partition $\{X_{p}\mid p\in P\}$
of a dense subset $X$ of $W$ satisfying $X_{p}^{\prime}\cap X=\bigcup_{q<p}X_{q}$
for $p\in P$. Here, $L\subseteq P$ defines which partition elements
are relatively compact, and the function $f$ specifies the size of
the finite partition elements.

We will consider the more general situation of a PO system $P$ and
a \emph{``trim $P$-partition''} $\{X_{p}\mid p\in P\}$ of a dense
subset of $W$. Trim partitions provide a physical representation
within the Stone space of the \emph{structure diagrams} defined by
Hanf~\cite{Hanf}; in particular, an $\omega$-Stone space is primitive
iff it admits a trim partition. We define the Cantor-Bendixson sequence
for a PO system $P$, and associated invariants $\nu(P)$ and a rank
function $r_{P}\colon K(P)\rightarrow\omega_{1}$, where $K(P)$ is
the perfect kernel of $P$, and show that $\nu(W)=\nu(P)$ and $r_{W}(x)=r_{P}(p)$
for $x\in X_{p}$ for any (primitive) $\omega$-Stone space~$W$
admitting a trim $P$-partition, so that $P$ significantly determines
the structure of~$W$. Moreover, if $P$ is a ``reduced'' PO system,
then we can determine $n(W)$ and $\rho(W)$ from $f$ and $L$ respectively.
A compact primitive $W$ is uniquely determined by the tuple $(P,f)$;
we show by example however that an extended PO system $(P,L,f)$ may
not uniquely determine a locally compact $\omega$-Stone space.

All scattered spaces are primitive, and so it is natural to ask how
to determine whether a uniform space is primitive. It is enough to
consider strongly uniform spaces $W$, which are uniquely determined
by their rank function $r_{W}$, which corresponds in turn to an additive
measure $\sigma\colon K(W)\rightarrow\omega_{1}$. The PI and primitivity
definitions readily extend to the notions of ``$\sigma$-PI'' and
``$\sigma$-primitive''. We show that a strongly uniform space $W$
is primitive iff $K(W)$ is $\sigma$-primitive, and further that
a compact open subset of $K(W)$ is $\sigma$-PI iff it is ``self-similar''
with respect to $\sigma$ at one of its points.

We use these ideas in the final section to develop a method for constructing
non-primitive spaces, by starting with a set of ``incompatible measures''
and constructing from them a measure that is not sufficiently self-similar,
and conclude with some specific examples. At the end we suggest several
potential areas for further study.

\textbf{Notation}: If $R$ is a Boolean ring and $A\in R$, we write
$(A)$ for $\{C\in R\mid C\subseteq A\}$, the ideal of $R$ generated
by $A$; $0$ for the zero element; and $1_{R}$ for the multiplicative
identity if $R$ is a Boolean algebra. We write $\mathscr{D}_{1}$
for the Cantor set and $\mathscr{D}_{0}$ for the (non-compact) Cantor
set minus a point. As usual, $\omega_{1}$ will denote the set of
all countable ordinals. If $W$ is a topological space and $S\subseteq W$,
we write $S^{\prime}$ for the \emph{derived set }of $S$, namely
the set of all limit points of $S$, and $\overline{S}$ for the closure
of $S$, so that $\overline{S}=S\cup S^{\prime}$.

We will use the following version of Vaught's Theorem twice, as cited
in~\cite[1.1.3]{PierceMonk}: 
\begin{theorem}[Vaught]
\label{(Vaught)}Suppose $\sim$ is a relation between the elements
of countable Boolean algebras $R$ and $S$ such that:

(i) $1_{R}\sim1_{S}$;

(ii) if $C\sim0_{S}$, then $C=0_{R}$; and vice versa;

(iii) if $C\sim D_{1}\dotplus D_{2}$ ($D_{1},D_{2}\in S$), then
we can write $C=C_{1}\dotplus C_{2}$, with $C_{i}\in R$ and $C_{i}\sim D_{i}$
$(i=1,2)$; and vice versa.

Then there is an isomorphism $\alpha\colon R\rightarrow S$ such that
each $C\in R$ can be expressed as $C=C_{1}\dotplus\ldots\dotplus C_{n}$
where $C_{i}\sim C_{i}\alpha$ for all $i\leqslant n$.
\end{theorem}

\section{\label{sec:Boolean-algebra-invariants}Invariants for $\omega$-Stone
spaces}

We recall the definitions of various invariants of compact $\omega$-Stone
spaces (e.g.\ see~\cite[1.5]{PierceMonk}) and extend them to general
$\omega$-Stone spaces. The invariants translate immediately to invariants
of the associated Boolean rings.
\begin{definition}
If $W$ is a topological space, its \emph{Cantor-Bendixson sequence}
is the list

\[
(W^{(0)},W^{(1)},\ldots,W^{(\xi)},\ldots),
\]

where $W^{(0)}=W$, $W^{(\xi+1)}=(W^{(\xi)})^{\prime}$, and $W^{(\eta)}=\bigcap_{\xi<\eta}W^{(\xi)}$
if $\eta$ is a limit ordinal.

If $W$ is an $\omega$-Stone space, then let:

\begin{eqnarray*}
\nu(W) & = & \min\{\xi\mid W^{(\xi)}=W^{(\xi+1)}\};\\
K(W) & = & W^{(\nu(W))}\text{, the \emph{perfect} \emph{kernel} of }W;\\
\lambda(W) & = & \min\{\xi\mid W^{(\xi)}-K(W)\text{ is closed}\};\\
n(W) & = & |W^{(\mu)}-K(W)|\text{ if }\nu(W)=\mu+1\text{ and }|W^{(\mu)}-K(W)|<\infty,\\
 & = & -\infty\text{ otherwise}.
\end{eqnarray*}

The \emph{rank function} of $W$ is the mapping $r_{W}\colon K(W)\rightarrow\omega_{1}$
given by $r_{W}(x)=\min\{\xi\in\omega_{1}\mid x\notin\overline{W^{(\xi)}-K(W)}\}$
for $x\in K(W)$.

$W$ is \emph{scattered }if $K(W)=\emptyset$ and is \emph{uniform}
if $\nu(W)=\lambda(W)$.

A map $f\colon W\rightarrow\omega_{1}$ is \emph{upper semi-continuous}
if $\{w\in W\mid f(w)\geqslant\xi\}$ is closed for all $\xi\in\omega_{1}$.
\end{definition}
\begin{remark}
\label{rem:invariant defns}$W^{\prime}$ is closed, so by transfinite
induction each $W^{(\xi)}$ is closed in $W$.

If $W$ is an $\omega$-Stone space and $K(W)$ is non-empty, then
$K(W)$ has no isolated points and so is homeomorphic to either $\mathscr{D}_{1}$
or $\mathscr{D}_{0}$.

For compact $W$, the usual definition is that $n(W)=|W^{(\mu)}-K(W)|$
if $\lambda(W)<\nu(W)$ and $\nu(W)=\mu+1$, and is $-\infty$ otherwise.
For non-compact $W$, however, there is the additional possibility
that $\nu(W)=\mu+1>\lambda(W)$ but $W^{(\mu)}-K(W)$ is infinite.
The revised definition caters for both these scenarios, as if $\nu(W)=\mu+1=\lambda(W)$
then $W^{(\mu)}-K(W)$ is not closed and so is infinite.
\end{remark}
\textbf{Notation}: for an ordinal $\xi$, let $W^{[\xi]}=W^{(\xi)}-W^{(\xi+1)}$
and let $K_{\xi}(W)=\{w\in K(W)\mid r_{W}(w)>\xi\}$, so that $r_{W}(w)=\min\{\xi\mid w\notin K_{\xi}(W)\}$.
We observe that $K_{\xi}(W)=K(W)\cap\overline{W^{(\xi)}-K(W)}$, and
is therefore closed in $W$, as $r_{W}(w)>\xi$ iff $w\in\overline{W^{(\xi)}-K(W)}$.

Parts~\ref{enu:prop1},~\ref{enu:prop4} and~\ref{enu:prop5} of
the next Proposition are simple extensions of~\cite[Propositions~1.4.1(e),~1.9.2]{PierceMonk}
to non-compact sets; we include proofs for completeness. 
\begin{proposition}
\label{Prop: invariant basics}\label{Lemma: Wx-Wx+1}Let $W$ be
an $\omega$-Stone space. Then:
\begin{enumerate}
\item \label{enu:prop1}if $A$ is a clopen subset of $W$, then $A^{(\xi)}=A\cap W^{(\xi)}$,
$K(A)=A\cap K(W)$ and $r_{A}(x)=r_{W}(x)$ for $x\in K(A)$;
\item \label{enu:prop2}$W^{[\xi]}{}^{\prime}-K(W)=W^{(\xi+1)}-K(W)$;
\item \label{enu:prop3}$\overline{W^{[\xi]}}=\overline{W^{(\xi)}-K(W)}=K_{\xi}(W)\cup(W^{(\xi)}-K(W))$;
\item \label{enu:prop4}$r_{W}$ is upper semi-continuous;
\item \label{enu:prop5}$\lambda(W)=\sup\{r_{W}(x)\mid x\in K(W)\}$.
\end{enumerate}
\end{proposition}
\begin{svmultproof}
For~\ref{enu:prop1}, we have $A^{(1)}=A^{\prime}=A\cap W^{(1)}$,
and $A\cap W^{(\xi)}$ is clopen in $W^{(\xi)}$. Hence $A^{(\xi)}=A\cap W^{(\xi)}$,
by transfinite induction, and so also $K(A)=A^{(\nu(W))}=A\cap K(W)$.
For the final statement, we observe that $A\cap\overline{W^{(\xi)}-K(W)}=\overline{A^{(\xi)}-K(A)}$.

For~\ref{enu:prop2}, clearly $W^{[\xi]}{}^{\prime}-K(W)\subseteq W^{(\xi+1)}-K(W)$.
If however $w\in W^{(\xi+1)}$ and $w\notin W^{[\xi]}{}^{\prime}$,
then we can find a compact open $A$ such that $w\in A$ and $A\cap W^{[\xi]}=\emptyset$.
So $A^{(\xi)}=A\cap W^{(\xi)}=A\cap W^{(\xi+1)}=A^{(\xi+1)}$ (using~\ref{enu:prop1}),
whence $w\in A^{(\xi+1)}=K(A)\subseteq K(W)$, as required.

For~\ref{enu:prop3}, it follows from~\ref{enu:prop2} that $\overline{W^{[\xi]}}-K(W)=W^{(\xi)}-K(W)$.
So $\overline{W^{[\xi]}}\subseteq\overline{W^{(\xi)}-K(W)}\subseteq\overline{W^{[\xi]}}$
and equality follows. Hence $\overline{W^{[\xi]}}=(\overline{W^{[\xi]}}\cap K(W))\cup(\overline{W^{[\xi]}}-K(W))=K_{\xi}(W)\cup(W^{(\xi)}-K(W))$.

For~\ref{enu:prop4}, we observed above that $K_{\xi}(W)$ is closed.
But $\{w\in K(W)\mid r_{W}(w)\geqslant\xi\}=\bigcap_{\eta<\xi}K_{\eta}(W)$,
and the result follows.

For~\ref{enu:prop5}, if $\xi=\lambda(W)$ then $K_{\xi}(W)=\emptyset$
and so $r_{W}(w)\leqslant\lambda(W)$ for all $w\in K(W)$. If however
$\xi<\lambda(W)$ then $K_{\xi}(W)\neq\emptyset$, as $W^{(\xi)}$
is closed, so we can find $w\in K(W)$ with $r_{W}(w)>\xi$.
\end{svmultproof}

\begin{corollary}
\label{Cor: alt defns}Let $W$ be an $\omega$-Stone space. Then
its invariants can be derived from $\{W^{[\xi]}\}$ as follows:

\begin{eqnarray*}
\nu(W) & = & \min\{\xi\mid W^{[\xi]}=\emptyset\};\\
\lambda(W) & = & \min\{\xi\mid\overline{W^{[\xi]}}\cap K(W)=\emptyset\}=\min\{\xi\mid K_{\xi}(W)=\emptyset\};\\
n(W) & = & |W^{[\mu]}|\text{ if }\nu(W)=\mu+1\text{ and }|W^{[\mu]}|<\infty,\\
 & = & -\infty\text{ otherwise};\\
r_{W}(x) & = & \min\{\xi\mid x\notin\overline{W^{[\xi]}}\}\text{ for }x\in K(W).
\end{eqnarray*}
\end{corollary}
\begin{svmultproof}
The formulae for $\nu(W)$ and $n(W)$ are immediate. That for $r_{W}(x)$
follows from Proposition~\ref{Prop: invariant basics}\ref{enu:prop3}.
Moreover, $W^{(\xi)}-K(W)$ is closed iff $\overline{W^{(\xi)}-K(W)}\cap K(W)=\overline{W^{[\xi]}}\cap K(W)=\emptyset$,
and the formula for $\lambda(W)$ follows. 
\end{svmultproof}

\begin{figure}
\centering{}\includegraphics[bb=200bp 0bp 900bp 540bp,scale=0.4]{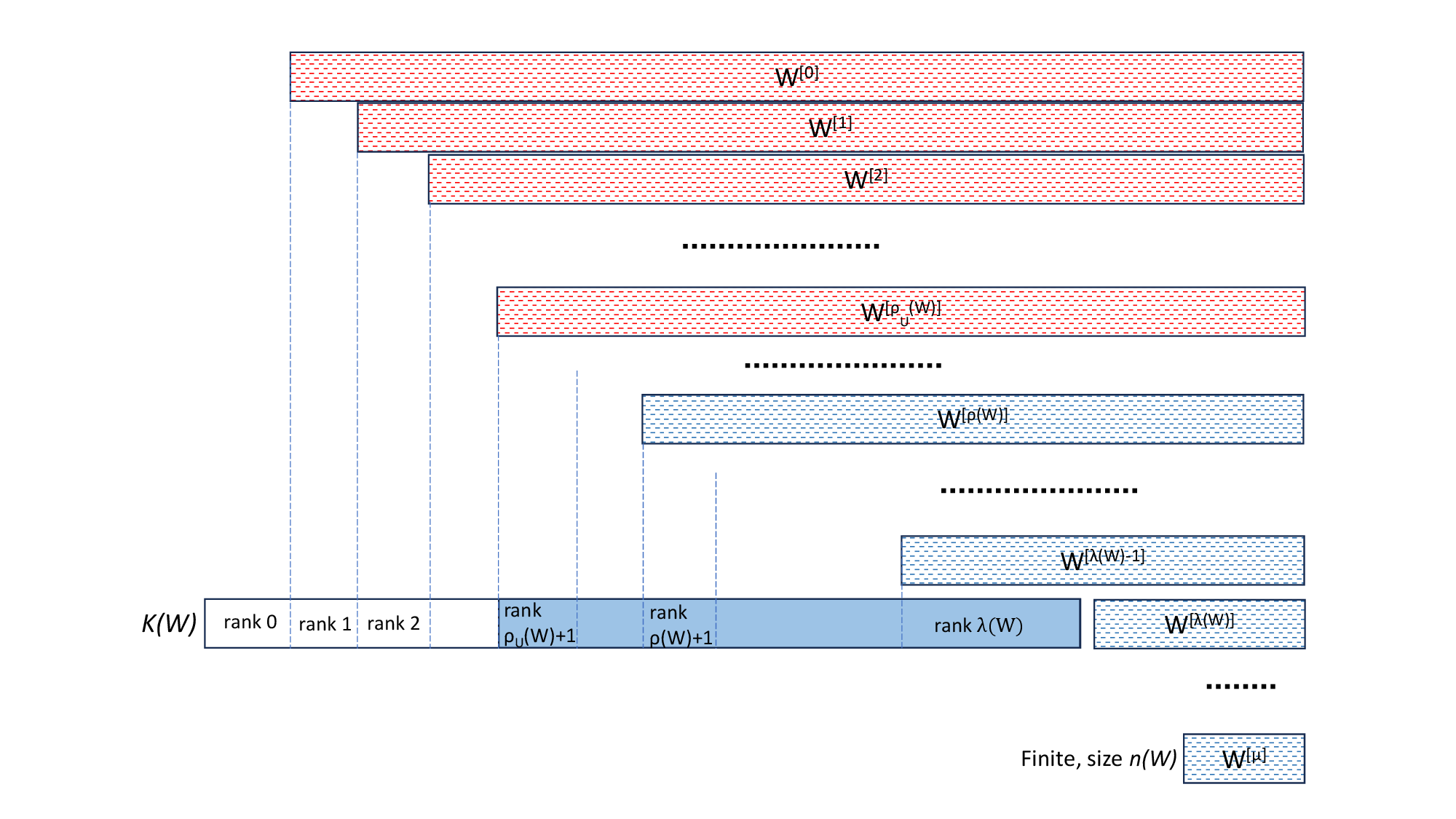}\caption{\label{fig:Structure}Structure of a typical Boolean space $W$}
\end{figure}

We must now add two further invariants to those described above, to
identify the compactness structure of $W$:
\begin{definition}
We define the \emph{compact rank }$\rho(W)$ and the \emph{uniform
compact rank }$\rho_{U}(W)$ of a general $\omega$-Stone space $W$
as follows: 
\begin{eqnarray*}
\rho(W) & = & \min\{\xi\mid\overline{W^{[\xi]}}\text{ is compact}\};\\
\rho_{U}(W) & = & \min\{\xi\mid K_{\xi}(W)\text{ is compact}\}.
\end{eqnarray*}
\end{definition}
\begin{remark}
It follows from Proposition~\ref{Prop: invariant basics} that $\rho_{U}(W)\leqslant\rho(W)\leqslant\nu(W)$
and that $\rho_{U}(W)\leqslant\lambda(W)$. Moreover, if $W$ is scattered
then $\rho(W)=\min\{\xi\mid W^{(\xi)}\text{ is compact}\}$, consistent
with the definition of compact rank in~\cite{FlumMartinez}.

Using Remark~\ref{rem:invariant defns}, we see that $n(W)=-\infty$
precisely when either $\lambda(W)=\nu(W)$ or $\rho(W)=\nu(W)$.

For any decomposition $W=Y\dotplus Z$, with $Y$ uniform and non-empty
and $Z$ scattered or empty, we will have $\rho_{U}(Y)=\rho_{U}(W)\leqslant\rho(Y)\leqslant\rho(W)$.
We will see in Section~\ref{sec:Structure-results-for} that we can
choose $Y$ to be ``strongly uniform'', with $\rho_{U}(W)=\rho_{U}(Y)=\rho(Y)$.
Thus $\rho_{U}(W)$ is the smallest possible value of $\rho(Y)$ for
any such uniform $Y$.
\end{remark}
Figure~\ref{fig:Structure} shows the structure of a typical Boolean
space for the case where $\rho_{U}(W)<\rho(W)<\lambda(W)<\nu(W)=\mu+1$
and $n(W)$ is finite. The discrete subspaces $W^{[\xi]}$ are shown
cross-hatched, and the closure of $W^{[\xi]}$ contains all subspaces
vertically below it. No attempt is made however to show the (potentially
complex) closure relations between the different rank subsets of $K(W)$.
Relatively compact subspaces are coloured blue and the others red. 

\section{\label{sec:Scattered--Stone-spaces}Scattered $\omega$-Stone spaces}

As a preliminary to the study of general $\omega$-Stone spaces, we
recap some existing results regarding \textbf{scattered }$\omega$-Stone
spaces. Briefly, each scattered $\omega$-Stone space $W$ is homeomorphic
to an ordinal space, and there is an elegant classification of these
in terms of the invariants $\nu(W)$, $\rho(W)$ and $n(W)$.

For non-zero ordinals $\alpha$ and $\mu$, let $T(\alpha)$ denote
the topological (Boolean) space consisting of $\alpha$ with the order
topology, and write $U(\mu)=[T(\omega^{\mu})]$ and $V(\mu)=[T(\omega^{\mu}+1)]$,
where $[W]$ denotes the homeomorphism class of a space $W$; also,
let $V(0)=[T(1)]$. Let the smallest term in the Cantor normal form
of a non-zero ordinal $\beta$ be $\omega^{h(\beta)}.m$ for some
$m\in\mathbb{N}_{+}$ ($h(\beta)$ is the \emph{Cantor Bendixson rank
}of $\beta$). We have the following (see e.g.\ \cite[Proposition~2.3.3]{Hilton}),
from which it follows immediately that all ordinal spaces are scattered:

\begin{equation}
T(\alpha)^{(\xi)}=\{\beta<\alpha\mid h(\beta)\geqslant\xi\}\label{ordCB}
\end{equation}

If $W=T(\alpha)$, it follows from~(\ref{ordCB}) that $W^{[\xi]}=\{\beta<\alpha\mid h(\beta)=\xi\}$,
so that (for example) $W^{[0]}$ includes $0$ and all successor ordinals
in $W$.

Let $[X].n$ denote the homeomorphism class of the disjoint union
of $n$ copies of a space $X$. In the literature, a compact ordinal
space is typically expressed as $T(\omega^{\mu}.n+1)$, which is homeomorphic
to $T(\omega^{\mu}+1).n$, with $T(\omega^{\mu}.n+1)^{(\mu)}=\{\omega^{\mu},\omega^{\mu}.2,\ldots,\omega^{\mu}.n\}$. 

The following result is due to Flum and Martinez~\cite[Theorem~2.2]{FlumMartinez};
we state it in the form adopted by Hilton~\cite[Theorem~2.4.1]{Hilton}:
\begin{theorem}[Flum, Martinez]
\label{FM}Every non-zero ordinal space $T(\alpha)$ belongs to one
and only one of the following homeomorphism classes, for some ordinals
$\mu$ and $\rho$, and $n\in\mathbb{N}_{+}$:
\begin{enumerate}
\item $V(\mu).n$ if $T(\alpha)$ is compact, where $\mu\geqslant0$;
\item $U(\mu)$ if none of the non-empty $T(\alpha)^{(\xi)}$ are finite
(equivalently, are compact), where $\mu>0$;
\item $V(\mu).n\dotplus U(\rho)$ if $T(\alpha)^{(\xi)}$ is compact precisely
when $\xi\geqslant\rho$ (i.e.\ $T(\alpha)$ has compact rank $\rho$),
where $0<\rho\leqslant\mu$.
\end{enumerate}
\end{theorem}
The uniqueness of $\mu$ for a given $T(\alpha)$, and of $n$ and
$\rho$ where relevant, follows by considering $|T(\alpha)^{(\xi)}|$
for varying ordinals $\xi$.

In the following definition, we note that if $W$ is scattered and
$\nu(W)$ is a limit ordinal then there is no minimum $W^{(\xi)}$
for $\xi<\nu(W)$, and none of these $W^{(\xi)}$ are compact as their
intersection is empty, and so $n(W)=-\infty$ and $\rho(W)=\nu(W)$.
If however $\nu(W)=\mu+1$, then $n(W)$ is finite or equal to $-\infty$
depending on whether $W^{(\mu)}$ is compact or non-compact (when
$\rho(W)=\nu(W)$). 
\begin{definition}
A scattered $\omega$-Stone space $W$ is of \emph{type }$(\nu(W),\nu(W))$
if $\rho(W)=\nu(W)$, otherwise $W$ is of \emph{type} $(\nu(W),\rho(W),n(W))$.
We write $\tau(W)$ for the type of $W$. 
\end{definition}
\begin{remark}
The empty set therefore has type $(0,0)$.
\end{remark}
\begin{proposition}
\label{Ordtable}For an ordinal $\alpha$, let $g(\alpha)$ denote
the smallest ordinal $\beta$ such that $T(\beta)\cong T(\alpha)$,
and write $\tau(\alpha)$ for $\tau(T(\alpha))$. Then $[T(\alpha)]$,
$g(\alpha)$ and $\tau(\alpha)$ are given by the following table,
where $\alpha=\omega^{\mu_{1}}.n_{1}+\cdots+\omega^{\mu_{k}}.n_{k}$
($\mu_{1}>\cdots>\mu_{k}$, each $n_{j}>0$) in Cantor normal form:

\begin{tabular}{>{\centering}p{3.18cm}ccc}
\hline 
$\alpha$  & $[T(\alpha)]$ & $g(\alpha)$ & $\tau(\alpha)$\tabularnewline
\hline 
\hline 
$0$ & $\emptyset$ & $0$ & $(0,0)$\tabularnewline
\hline 
$n$\\
($n\geqslant1)$ & $V(0).n$ & $n$ & $(1,0,n)$\tabularnewline
\hline 
$\omega^{\mu}$\\
($\mu>0)$ & $U(\mu)$ & $\omega^{\mu}$ & $(\mu,\mu)$\tabularnewline
\hline 
$\omega^{\mu}.(n+1)$ ($\mu>0,n\geqslant1$) & $V(\mu).n\dotplus U(\mu)$ & $\omega^{\mu}.(n+1)$ & $(\mu+1,\mu,n)$\tabularnewline
\hline 
$\omega^{\mu_{1}}.n_{1}+\cdots+1.n_{k}$ ($k>1$, $\mu_{k}=0$) & $V(\mu_{1}).n_{1}$ & $\omega^{\mu_{1}}.n_{1}+1$ & $(\mu_{1}+1,0,n_{1})$\tabularnewline
\hline 
$\omega^{\mu_{1}}.n_{1}+\cdots+\omega^{\mu_{k}}.n_{k}$ ($k>1,\mu_{k}>0$) & $V(\mu_{1}).n_{1}\dotplus U(\mu_{k})$ & $\omega^{\mu_{1}}.n_{1}+\omega^{\mu_{k}}$ & $(\mu_{1}+1,\mu_{k},n_{1})$\tabularnewline
\hline 
\end{tabular}
\end{proposition}
\begin{svmultproof}
The formulae for $[T(\alpha)]$ follow immediately from Flum and Martinez~\cite[p.~790]{FlumMartinez},
and those for $g(\alpha)$ follow from those formulae by inspection. 

Now $T(\omega^{\mu})^{(\xi)}=\bigcup\{T(\omega^{\eta}.n+1)^{(\xi)}\mid n\in\mathbb{N}_{+},\xi\leqslant\eta<\mu\}$,
and $T(\alpha)^{(\xi)}=T(\alpha)\cap T(\omega^{\mu})^{(\xi)}$ which
is an open subset of $T(\omega^{\mu})^{(\xi)}$ for $\alpha<\omega^{\mu}$.
But for fixed $\xi<\mu$, $\{T(\omega^{\eta}.n+1)^{(\xi)}\mid n\in\mathbb{N}_{+},\xi\leqslant\eta<\mu\}$
are non-empty and strictly increasing, as $T(\omega^{\eta}.n+1)^{(\xi)}$
contains $\omega^{\eta}.n$ but not $\omega^{\eta}.(n+1)$, by~(\ref{ordCB}).
Hence none of the $T(\omega^{\mu})^{(\xi)}$ are compact for $\xi<\mu$,
and so $\omega^{\mu}$ has type $(\mu,\mu)$ as $T(\omega^{\mu})^{(\mu)}=\emptyset$. 

Furthermore, $T(\omega^{\mu}.n+1)^{(\mu)}=\{\omega^{\mu},\omega^{\mu}.2,\ldots,\omega^{\mu}.n\}$
by~(\ref{ordCB}), and so $T(\omega^{\mu}.n+1)$ has type $(\mu+1,0,n)$.
Using these two results together with the fact that $(W\dotplus X)^{(\xi)}=W^{(\xi)}\dotplus X^{(\xi)}$,
we see that the type of $V(\mu).n\dotplus U(\rho)$ is $(\mu+1,\rho,n)$
for $0<\rho\leqslant\mu$.
\end{svmultproof}

\begin{remark}
If we write $U(0)=V(\mu).0=\emptyset$ and $\omega^{0}=1$, then we
have $[T(\omega^{\mu_{1}}.n_{1}+\cdots+\omega^{\mu_{k}}.n_{k})]=V(\mu_{1}).(n_{1}-e(k,\mu_{1}))\dotplus U(\mu_{k})$,
where $e(k,\mu_{1})=1$ if $k=1$ and $\mu_{1}>0$, and otherwise
$e(k,\mu_{1})=0$.
\end{remark}
\begin{corollary}[Flum, Martinez]
\label{Ordspacetypes}An ordinal space is determined up to homeomorphism
by its type. Moreover, every possible type $(\mu,\mu)$ and $(\mu+1,\rho,n)$
arises for $0\leqslant\rho\leqslant\mu\in\omega_{1}$ and $n\in\mathbb{N}_{+}$.
\end{corollary}
\begin{svmultproof}
Immediate from Proposition~\ref{Ordtable}.
\end{svmultproof}

\begin{corollary}
The (commutative) monoid $S$ of homeomorphism classes of ordinal
spaces under disjoint union is generated by $\{V(0),U(\mu),V(\mu)\mid\mu>0\}$,
with the relations:

\begin{align*}
U(\mu)\dotplus U(\nu) & =U(\max(\mu,\nu));\\
V(\mu)\dotplus U(\nu) & =U(\nu)\text{ if }\mu<\nu;\\
V(\mu)\dotplus V(\nu) & =V(\nu)\text{ if }\mu<\nu.
\end{align*}

Equivalently, $S=\{(\mu,\mu),(\mu+1,\rho,n)\mid0\leqslant\rho\leqslant\mu,n\in\mathbb{N}_{+}\}$,
with the following addition rules:

\begin{align*}
(\mu_{1},\mu_{1})+(\mu_{2},\mu_{2}) & =(\max(\mu_{1},\mu_{2}),\max(\mu_{1},\mu_{2}));\\
(\mu+1,\rho_{1},n_{1})+(\mu+1,\rho_{2},n_{2}) & =(\mu+1,\max(\rho_{1},\rho_{2}),n_{1}+n_{2});\\
(\mu_{1}+1,\rho_{1},n_{1})+(\mu_{2}+1,\rho_{2},n_{2}) & =(\mu_{2}+1,\max(\rho_{1},\rho_{2}),n_{2})\text{ if }\mu_{1}<\mu_{2};\\
(\mu_{1}+1,\rho_{1},n_{1})+(\mu_{2},\mu_{2}) & =(\mu_{1}+1,\max(\rho_{1},\mu_{2}),n_{1})\text{ if }\mu_{1}\geqslant\mu_{2};\\
(\mu_{1}+1,\rho_{1},n_{1})+(\mu_{2},\mu_{2}) & =(\mu_{2},\mu_{2})\text{ if }\mu_{1}<\mu_{2}.
\end{align*}
\end{corollary}
\begin{svmultproof}
Let $W$ and $X$ be ordinal spaces of type $U(\mu)$ or $V(\nu)$.
The first set of formulae are easily derived by considering the cardinality
and compactness of each $W^{(\xi)}\dotplus X^{(\xi)}$ and applying
Theorem~\ref{FM} and Proposition~\ref{Ordtable}. The remaining
statements then follow from Proposition~\ref{Ordtable} and the first
set of formulae. 
\end{svmultproof}

Flum and Martinez also proved that any countable scattered space is
homeomorphic to an ordinal space. We provide an elementary proof of
this using Vaught's Theorem~\ref{(Vaught)}; the corresponding proofs
in~\cite[Theorem~2.2]{FlumMartinez} for the general case and in~\cite{PierceMonk}
for the compact case use model-theoretic and ring-theoretic uniqueness
results respectively. 
\begin{theorem}[{Flum and Martinez~\cite[Theorem~2.3]{FlumMartinez}}]
\label{Thm: Flum Martinez}The scattered $\omega$-Stone space~$W$
is homeomorphic to the (unique) ordinal space of type $\tau(W)$.
\end{theorem}
\begin{svmultproof}
It is enough to show that two $\omega$-Stone spaces of the same type
are homeomorphic, as every possible type is realised by an ordinal
space, by Corollary~\ref{Ordspacetypes}. So let $W$ and $X$ be
$\omega$-Stone spaces such that $\tau(W)=\tau(X)$, and let $R$
and $S$ denote the compact open subsets of $W$ and $X$ respectively. 

Suppose first that $W$ is compact, so that $\rho(X)=\rho(W)=0$ and
$X$ is also compact. For $A\in R$ and $B\in S$, let $A\sim B$
if $\tau(A)=\tau(B)$. We must show that $\sim$ satisfies the Vaught
conditions of Theorem~\ref{(Vaught)}. Clearly $W=1_{R}\sim1_{S}=X$,
and $A\sim0$ iff $A=0$. Suppose then that $C\sim D=D_{1}\dotplus D_{2}$,
where $D,D_{i}\in S$ for $i=1,2$ and $\nu(D_{1})\leqslant\nu(D_{2})$:
say $\tau(C)=(\mu+1,0,n)$ and $\nu(D_{1})=\mu_{1}+1$ with $\mu_{1}\leqslant\mu$.
If $\mu_{1}=\mu$, then $\nu(D_{1})=\nu(D_{2})=\mu$ and so $n=n(D_{1})+n(D_{2})$.

We can choose $C_{1}\subseteq C-C^{(\mu_{1}+1)}$ such that $C_{1}\in R$
and $|C_{1}\cap C^{[\mu_{1}]}|=n(D_{1})$: as $C^{(\mu_{1}+1)}$ is
closed, and $C^{[\mu_{1}]}$ is discrete infinite if $\mu_{1}<\mu$
and of size $n$ if $\mu_{1}=\mu$. So $\nu(C_{1})=\mu_{1}+1$, as
$C_{1}^{[\mu_{1}]}=C_{1}\cap C^{[\mu_{1}]}\neq\emptyset$ and $C_{1}^{(\mu_{1}+1)}=\emptyset$,
whence $\tau(C_{1})=\tau(D_{1})$. 

Let $C_{2}=C-C_{1}$. Then $\tau(C_{2})=\tau(C)=\tau(D)=\tau(D_{2})$
if $\mu_{1}<\mu$, while $\tau(C_{2})=(\mu+1,0,n-n(D_{1}))=\tau(D_{2})$
if $\mu_{1}=\mu$. So condition (iii) is satisfied, and we have $W\cong X$.

Suppose now that $W$ is not compact. Choose $A\in R$ and $B\in S$
such that $W^{(\rho)}\subseteq A$ and $X^{(\rho)}\subseteq B$, where
$\rho=\rho(W)>0$; if $\rho=\nu(W)$, let $A=B=\emptyset$. If $\rho<\nu(W)$,
with $\nu(W)=\mu+1$, say, then $W^{[\mu]}\subseteq A$, so $\tau(A)=(\nu(W),0,n(W))=\tau(B)$.
So we can find an isomorphism $\theta_{0}\colon(A)\rightarrow(B)$;
this is trivially true if $\rho=\nu(W)$.

Suppose now that we have $C\in R$ and $D\in S$ with $A\subseteq C$
and $B\subseteq D$, and that we have extended $\theta_{0}$ to an
isomorphism $\theta_{1}\colon(C)\rightarrow(D)$. Suppose also that
$E\in R$ such that $E\cap C=\emptyset$. Then $E^{(\rho)}=\emptyset$:
say $\tau(E)=(\xi+1,0,n)$, where $\xi+1\leqslant\rho$. Now $\overline{X^{[\xi]}}$
is not compact, so we can find $F\in S$ such that $F\cap D=\emptyset$,
$F\subseteq X-X^{(\xi+1)}$ and $|F\cap X^{[\xi]}|=n$. Then $\tau(F)=\tau(E)$,
so we can extend $\theta_{1}$ to an isomorphism $\theta_{2}\colon(C\dotplus E)\rightarrow(D\dotplus F)$. 

A routine back and forth argument now extends $\theta_{0}$ to an
isomorphism between $R$ and $S$, and so $W\cong X$, as required.
\end{svmultproof}

\section{\label{sec:Structure-results-for}Structure results for general $\omega$-Stone
spaces}

We now turn from scattered spaces to general spaces. For compact $\omega$-Stone
spaces, we recall the following decomposition result of Ketonen and
uniqueness and existence results of Pierce. The goal of this section
is to derive related decomposition, uniqueness and existence results
for locally compact $\omega$-Stone spaces.
\begin{theorem}[Ketonen~\cite{Ketonen}]
\label{(Ketonen)}If $X$ is a non-uniform compact $\omega$-Stone
space, then we can write $X=Y\dotplus Z$, with $Y$ and $Z$ clopen
in $X$, where:
\begin{enumerate}
\item $Z$ is scattered, with $\nu(Z)=\nu(X)$ and $n(Z)=n(X)$,
\item either $Y$ is empty, or $Y$ is uniform with $\lambda(Y)=\lambda(X)$
and $K(Y)=K(X)$.
\end{enumerate}
\end{theorem}
\begin{remark}
The decomposition is unique up to homeomorphism of $Y$ and $Z$.
\end{remark}
\textbf{Notation}: if $U$ and $V$ are topological spaces and $r\colon U\rightarrow\omega_{1}$
and $s\colon V\rightarrow\omega_{1}$ are maps, we will write $(U,r)\cong(V,s)$
if there is a homeomorphism $\phi\colon U\rightarrow V$ such that
$s(u\phi)=r(u)$ for all $u\in U$.
\begin{theorem}[{Pierce~\cite[Theorem~1.10.1]{PierceMonk}}]
\label{Thm: invariants uniqueness}If $X$ and $Y$ are compact $\omega$-Stone
spaces, then $X$ is homeomorphic to $Y$ iff $\nu(X)=\nu(Y)$, $\lambda(X)=\lambda(Y)$,
$n(X)=n(Y)$ and $(K(X),r_{X})\cong(K(Y),r_{Y})$. 
\end{theorem}
\begin{remark}
The proof in~\cite{PierceMonk} uses Theorem~\ref{(Ketonen)} together
with Vaught's Theorem~\ref{(Vaught)}. In particular, compact uniform
$\omega$-Stone spaces are homeomorphic iff they have isomorphic rank
functions. The requirement that $\lambda(X)=\lambda(Y)$ can be omitted,
as by Proposition~\ref{Prop: invariant basics}\ref{enu:prop5} the
value of $\ensuremath{\lambda(X)}$ can be deduced from the rank function
$r_{X}$.
\end{remark}
\begin{theorem}[{Pierce~\cite[Theorem~1.11.1,~Corollary~1.11.2]{PierceMonk}}]
\label{Thm: Pierce existence}
\begin{enumerate}
\item If $r\colon\mathscr{D}_{1}\rightarrow\omega_{1}$ is upper semi-continuous,
then there is a uniform $\omega$-Stone space $W$ such that $(K(W),r_{W})\cong(\mathscr{D}_{1},r)$;
\item if $n\in\mathbb{N}_{+}$, $\mu<\omega_{1}$ and $r\colon\mathscr{D}_{1}\rightarrow\omega_{1}$
is upper semi-continuous such that $r(x)\leqslant\mu$ for all $x\in\mathscr{D}_{1}$,
then there is an $\omega$-Stone space $W$ such that $\nu(W)=\mu+1$,
$(K(W),r_{W})\cong(\mathscr{D}_{1},r)$ and $n(W)=n$.
\end{enumerate}
\end{theorem}
To extend these results to a general non-compact $\omega$-Stone space,
we will need to factor in the compact ranks defined earlier. We will
also frequently use the following:
\begin{lemma}
\label{Lemma: dp-results}Let $W$ be an $\omega$-Stone space.
\begin{enumerate}
\item \label{enu:dp1}If $A\subseteq B\subseteq W$, with $A$ and $B$
closed subsets of $W$ and $A$ compact and open in $B$, then we
can find a compact open subset $C$ of $W$ such that $A=B\cap C$;
\item \label{enu:dp2}if $W$ is not compact, then we can write $W=B_{1}\dotplus B_{2}\dotplus\cdots$,
where the $B_{j}$ are disjoint compact open subsets of $W$;
\item \label{enu:dp3}if $W=\dotplus\{B_{j}\mid j\in J\}$, where the $B_{j}$
are disjoint compact open subsets of $W$and $J$ is finite or countably
infinite, then:
\end{enumerate}
\begin{eqnarray*}
W^{(\xi)} & = & \dotplus\{B_{j}^{(\xi)}\mid j\in J\};\\
W^{[\xi]} & = & \dotplus\{B_{j}^{[\xi]}\mid j\in J\};\\
\nu(W) & = & \sup\{\nu(B_{j})\mid j\in J\};\\
K(W) & = & \dotplus\{K(B_{j})\mid j\in J\};\\
K_{\xi}(W) & = & \dotplus\{K_{\xi}(B_{j})\mid j\in J\};\\
\lambda(W) & = & \sup\{\lambda(B_{j})\mid j\in J\};\\
\rho(W) & = & \max\{\rho(B_{j})\mid j\in J\}\text{\text{ if J is finite; and}}\\
\rho_{U}(W) & = & \max\{\rho_{U}(B_{j})\mid j\in J\}\text{ if J is finite.}
\end{eqnarray*}
\end{lemma}
\begin{svmultproof}
Let $R$ be the Boolean ring of compact open subsets of $W$.

\ref{enu:dp1}~We can write $A=B\cap D$, where $D$ is open in $W$.
Write $D$ as a (possibly infinite) union of elements of $R$. By
compactness, a finite subset of these elements will cover $A$.

\ref{enu:dp2}~Let $\{A_{j}\mid j\geqslant1\}$ be an enumeration
of $R$. Define $\{B_{n}\}$ iteratively by setting $B_{n}=\bigcup_{j\leqslant n}A_{j}-\bigcup_{j<n}B_{j}$,
and discarding any empty $B_{n}$.

\ref{enu:dp3}~It is easy to see that $W^{\prime}=\dotplus\{B_{j}^{\prime}\mid j\in J\}$,
from which the first 4 assertions readily follow. Moreover, $r_{B_{j}}(x)=r_{W}(x)$
for $x\in B_{j}$, from which we obtain the formula for $K_{\xi}(W)$.
Next, $\lambda(W)=\sup\{r_{W}(x)\mid x\in K(W)\}=\sup\{r_{B_{j}}(x)\mid x\in K(B_{j}),j\in J\}=\sup\{\lambda(B_{j})\mid j\in J\}$.
The final two assertions follow immediately from the previous assertions
and from the definitions of $\rho(W)$ and $\rho_{U}(W)$.
\end{svmultproof}

\begin{lemma}
\label{uniformaddition}If $W$ and $X$ are $\omega$-Stone spaces,
with $W$ uniform and $\nu(X)\leqslant\nu(W)$, then $W\dotplus X$
is uniform.
\end{lemma}
\begin{svmultproof}
$\nu(W\dotplus X)\geqslant\lambda(W\dotplus X)\geqslant\lambda(W)=\nu(W)=\nu(W\dotplus X)$.
\end{svmultproof}

The next definition and proposition provide the key bridge between
the compact and non-compact cases.
\begin{definition}
A uniform $\omega$-Stone space $W$ is \emph{strongly uniform} if
also $\rho_{U}(W)=\rho(W)$.
\end{definition}
\begin{proposition}
\label{Stronglyuniformlift}Let $W$ be a uniform $\omega$-Stone
space, and let $K(W)=\dotplus\{J_{n}\mid n<\alpha\}$ be a decomposition
of $K(W)$ into compact open subsets of $K(W)$, where $\alpha\leqslant\omega$.
Then there is a decomposition $W=\dotplus\{A_{n}\mid n<\alpha\}\dotplus Y$
of $W$ such that 
\begin{enumerate}
\item for each $n<\alpha$, $A_{n}$ is a uniform compact open subset with
$K(A_{n})=J_{n}$;
\item $Y$ is empty if $W$ is strongly uniform;
\item otherwise $Y$ is a scattered clopen subset of $W$ of type $(\rho(W),\rho(W))$.
\end{enumerate}
\end{proposition}
\begin{svmultproof}
Suppose first that $W$ is not compact, and let $W=B_{1}\dotplus B_{2}\dotplus\cdots$
be a decomposition of $W$ into disjoint compact open subsets. By
grouping together the $B_{j}$ and using the compactness of each $J_{n}$,
we may assume that $J_{n}\subseteq\bigcup_{1\leqslant j\leqslant n}B_{j}$
for $n<\alpha$. 

\textbf{Step 1}

We claim that we can find disjoint compact open subsets $\{C_{i}\mid i<\alpha\}$
of $W$ such that for all $i$: (i) $C_{i}\subseteq\bigcup_{1\leqslant j\leqslant i}B_{j}$,
(ii) $K(C_{i})=J_{i}$, (iii) if $J_{i}\cap B_{j}=\emptyset$ then
$C_{i}\cap B_{j}=\emptyset$ for $j\leqslant i$; (iv) $C_{i}$ is
uniform. 

For suppose $1\leqslant n+1<\alpha$ and we already have disjoint
compact open subsets $\{C_{1},\ldots,C_{n}\}$ of $W$ satisfying
(i) to (iv) for each $i\leqslant n$. Find a compact open subset $C_{n+1}$
of $W$ such that $C_{n+1}\cap K(W)=J_{n+1}$, so that $K(C_{n+1})=J_{n+1}$.
Replacing $C_{n+1}$ by

\[
C_{n+1}\cap\left\{ \bigcup_{j\leqslant n+1}\{B_{j}\mid B_{j}\cap J_{n+1}\neq\emptyset\}-\bigcup_{i\leqslant n}C_{i}\right\} 
\]

we obtain conditions (i) and (iii) while preserving $C_{n+1}\cap K(W)=J_{n+1}$,
as $C_{i}\cap J_{n+1}=\emptyset$ for $i\leqslant n$. Lastly, we
can use Theorem~\ref{(Ketonen)} to write $C_{n+1}$ as the union
of a uniform and scattered component, and discard the scattered component
to obtain that $C_{n+1}$ is uniform, as a scattered space has an
empty perfect kernel.

\textbf{Step 2}

For each $j$, $K(B_{j})$ is compact and meets only finitely many
of the $J_{n}$, and so $B_{j}$ meets only finitely many of the $C_{n}$
(using condition (iii) for $n\geqslant j$ and so for all $n$, as
$J_{n}\cap B_{j}=J_{n}\cap K(B_{j}))$. Let $D_{j}=B_{j}-\bigcup_{n<\alpha}C_{n}$,
which is compact and open, and is scattered as $K(W)\subseteq\bigcup_{n<\alpha}C_{n}$:
say $\nu(D_{j})=\eta_{j}+1$ for $D_{j}\neq\emptyset$, with $\eta_{j}<\nu(W)=\lambda(W)$. 

Let $\rho=\rho(W)$ and $\sigma=\rho_{U}(W)$. As $\overline{W^{[\rho]}}$
and $K_{\sigma}(W)$ are compact, we can find $N\in\mathbb{N}$ such
that $B_{j}\cap\overline{W^{[\rho]}}=\overline{B_{j}^{[\rho]}}=\emptyset$
and $C_{j}\cap K_{\sigma}(W)=K_{\sigma}(C_{j})=\emptyset$ for $j>N$.
Now if $D_{j}\neq\emptyset$ and $j>N$, then $|D_{j}^{[\eta_{j}]}|>0$,
and so $\eta_{j}<\rho$. 

\textbf{Step 3}

Let $L=\{j\in\mathbb{N}_{+}\mid D_{j}\neq\emptyset\wedge\sigma\leqslant\eta_{j}<\rho\}$
and $M=\{j\in\mathbb{N}_{+}-L\mid D_{j}\neq\emptyset\}$. We define
a function $f\colon M\rightarrow\mathbb{N}$ as follows, which will
enable the uniform $\{C_{n}\}$ to ``swallow up'' the scattered
$\{D_{j}\}$ for $j\in M$:

(a) if $\eta_{j}\geqslant\rho$ (so that $j\leqslant N$), let $f(j)$
be the smallest value of $n$ such that $K_{\eta_{j}}(W)\cap J_{n}\neq\emptyset$;

(b) if $\eta_{j}<\sigma$, let $f(j)$ be the smallest value of $n>j$
such that $K_{\eta_{j}}(W)\cap J_{n}\neq\emptyset$; this is well
defined, as $K_{\xi}(W)$ is not compact for $\xi<\sigma$.

We observe firstly that if $j\in M$ then $J_{f(j)}$ contains elements
with rank greater than $\eta_{j}$, and so $\lambda(C_{f(j)})\geqslant\nu(D_{j})$;
and secondly that $f^{-1}(n)$ is finite for each $n>0$, as if $f(j)=n$
then either $j\leqslant N$ or $j<n$.

\textbf{Step 4}

Finally, let $A_{n}=C_{n}\cup\bigcup_{j\in M}\{D_{j}\mid f(j)=n\}$
for $n<\alpha$, which is compact, open and uniform by Lemma~\ref{uniformaddition},
and let $Y=\dotplus\{D_{j}\mid j\in L\}$. Then $W=\dotplus\{A_{n}\mid n<\alpha\}\dotplus Y$
and $K(A_{n})=J_{n}$ for each $n<\alpha$. 

If $W$ is strongly uniform, then $L$ and $Y$ are empty. Otherwise,
if $\sigma\leqslant\xi<\rho$, then $K_{\sigma}(C_{n})=\emptyset$,
$\lambda(C_{n})\leqslant\sigma$ and $C_{n}^{[\xi]}=\emptyset$ for
$n>N$. So if $j>N$ then $D_{j}^{[\xi]}=B_{j}^{[\xi]}$, as $B_{j}\cap C_{n}=\emptyset$
for $n<j$ and so $D_{j}=B_{j}-\bigcup_{n\geqslant j}C_{n}$. But
$\overline{W^{[\xi]}}$ is not compact, so we can find infinitely
many $j>N$ such that $B_{j}\cap W^{[\xi]}\neq\emptyset$, whence
also $D_{j}^{[\xi]}\neq\emptyset$. Therefore $\rho(Y)>\xi$. But
$\nu(Y)\leqslant\rho$. It follows that $\rho(Y)=\nu(Y)=\rho$, and
$Y$ is scattered of type $(\rho,\rho)$, as required.

Finally, if $W$ is compact, so that $\alpha\in\mathbb{N}_{+}$ and
$W$ is strongly uniform, then a simpler version of step 1 yields
$W=\dotplus\{A_{n}\mid n<\alpha\}\dotplus Y$, where $A_{n}$ is uniform,
$K(A_{n})=J_{n}$ and $Y$ is compact and scattered. By Lemma~\ref{Lemma: dp-results},
$\nu(W)=\lambda(W)=\lambda(A_{n})=\nu(A_{n})$ for some $n<\alpha$,
and so by Lemma~\ref{uniformaddition} we can replace $A_{n}$ with
$A_{n}\dotplus Y$ to obtain that $Y$ is empty.
\end{svmultproof}

We will now apply the previous Proposition to obtain structure results
for increasingly general types of $\omega$-Stone spaces, starting
with strongly uniform spaces.
\begin{corollary}
\label{Cor:dpcompact uniform}An $\omega$-Stone space is strongly
uniform iff it is a direct product of compact open uniform spaces.
\end{corollary}
\begin{svmultproof}
``Only if'' is an immediate consequence of Proposition~\ref{Stronglyuniformlift}.
Suppose then that $W=A_{1}\dotplus A_{2}\dotplus\cdots$, with each
$A_{n}$ a compact open uniform subset of $W$. By Lemma~\ref{Lemma: dp-results},
$\lambda(W)=\sup_{n\geqslant1}\{\lambda(A_{n})\}=\sup_{n\geqslant1}\{\nu(A_{n})\}=\nu(W)$,
so that $W$ is uniform. We must show that $\rho_{U}(W)=\rho(W)$.
But if $\xi<\rho(W)$, then $W^{[\xi]}$ is not compact. So there
are infinitely many $n$ such that $A_{n}^{[\xi]}\neq\emptyset$,
and hence such that $\lambda(A_{n})=\nu(A_{n})>\xi$, so that $K_{\xi}(A_{n})\neq\emptyset$.
Therefore $K_{\xi}(W)$ is not compact and $\xi<\rho_{U}(W)$. Thus
$\rho_{U}(W)=\rho(W)$, as required.
\end{svmultproof}

\begin{theorem}
\label{Thm: strongly uniform unique}Let $W$ be a strongly uniform
$\omega$-Stone space. Then $W$ is determined up to homeomorphism
by the isomorphism type of the pair $(K(W),r_{W})$. Moreover, the
map $W\mapsto(K(W),r_{W})$ induces a bijection between the homeomorphism
classes of strongly uniform $\omega$-Stone spaces and isomorphism
classes of pairs $(K,r)$, where $K\in\{\mathscr{D}_{0},\mathscr{D}_{1}\}$
and $r\colon K\rightarrow\omega_{1}$ is upper semi-continuous.
\end{theorem}
\begin{svmultproof}
Suppose $W$ and $X$ are strongly uniform $\omega$-Stone spaces
and that there is an isomorphism $\theta\colon K(W)\rightarrow K(X)$
such that $r_{X}(w\theta)=r_{W}(w)$ for all $w\in K(W)$. Let $K(W)=J_{1}\dotplus J_{2}\dotplus\cdots$
be a decomposition of $K(W)$ into compact open subsets of $K(W)$,
so that $K(X)=J_{1}\theta\dotplus J_{2}\theta\dotplus\cdots$.

Apply Proposition~\ref{Stronglyuniformlift} to obtain decompositions
$W=A_{1}\dotplus A_{2}\dotplus\cdots$ and $X=B_{1}\dotplus B_{2}\dotplus\cdots$
into disjoint compact open uniform subsets such that $K(A_{n})=J_{n}$
and $K(B_{n})=J_{n}\theta$ for each $n$. Now $r_{B_{n}}(w\theta)=r_{A_{n}}(w)$
for all $w\in J_{n}$, so by Theorem~\ref{Thm: invariants uniqueness}
$A_{n}\cong B_{n}$ for each $n\geqslant1$, and we can combine these
homeomorphisms to obtain the desired homeomorphism from $W$ to $X$.

For the second statement, it is immediate that if $W\cong X$ then
$(K(W),r_{W})\cong(K(X),r_{X})$, so it remains to show that every
possible pair $(K,r)$ can arise. If $K$ is compact, then the result
follows from Theorem~\ref{Thm: Pierce existence}. Suppose instead
that we are given a pair $(K,r)$ with $K$ non-compact, and write
$K=K_{1}\dotplus K_{2}\dotplus\cdots$ with each $K_{n}$ being a
compact open subset of $K$. For each $n$ apply Theorem~\ref{Thm: Pierce existence}
to find a compact uniform space $A_{n}$ such that $(K_{n},r_{n})\cong(K(A_{n}),r_{A_{n}})$,
where $r_{n}$ is the restriction of $r$ to $K_{n}$. Let $W$ be
the disjoint union of the $A_{n}$, which is strongly uniform by Corollary~\ref{Cor:dpcompact uniform},
with $(K(W),r_{W})\cong(K,r)$ as required. 
\end{svmultproof}

\begin{remark}
For any such decomposition $W=A_{1}\dotplus A_{2}\dotplus\cdots$
into compact open uniform spaces, there will only be finitely many
$A_{n}$ such that $\nu(A_{n})>\rho(W)$. Hence we can write $W=X\dotplus Z$,
where $X$ is strongly uniform with $\rho(X)=\nu(X)=\rho(W)$ and
$Z$ is compact and uniform with $\nu(Z)=\nu(W)$. For primitive spaces
(see Section~\ref{sec:Primitive-spaces,-trim} below), it can be
shown that there is a ``minimal'' such $Z$, unique up to homeomorphism,
but $X$ may not be unique: e.g.\ see the space~$Y$ in Example~\ref{exa: non-unique PLf}.
For non-primitive spaces, a variation of Example~\ref{exa:incompat}
shows that the homeomorphism class of $Z$ may not be unique. 
\end{remark}
\begin{theorem}
\label{Thm: uniform decomposn}An $\omega$-Stone space $W$ is uniform
iff we can write $W=X\dotplus Y$, where $X$ is strongly uniform
and $Y$ is either empty or is a scattered space of type $(\rho,\rho)$,
for some $\rho_{U}(X)<\rho\leqslant\nu(X)$ (and where $\rho=\rho(W)$).
\end{theorem}
\begin{svmultproof}
The ``if'' statement follows immediately from Lemma~\ref{uniformaddition}.

Suppose instead that $W$ is uniform, and use Proposition~\ref{Stronglyuniformlift}
to write $W=Y\dotplus A_{1}\dotplus A_{2}\dotplus\cdots$, with each
$A_{n}$ a compact open subset of $W$ and $Y$ either empty (if $W$
is strongly uniform) or scattered of type $(\rho(W),\rho(W))$. Let
$X=A_{1}\dotplus A_{2}\dotplus\cdots$, which is strongly uniform
by Corollary~\ref{Cor:dpcompact uniform}, with $\nu(X)=\lambda(X)=\lambda(W)=\nu(W)$.
If now $Y$ is non-empty, then $\rho_{U}(X)=\rho_{U}(W)<\rho(W)\leqslant\nu(W)=\nu(X)$,
as required.
\end{svmultproof}

\textbf{Notation}: if $K\in\{\mathscr{D}_{0},\mathscr{D}_{1}\}$ and
$r\colon K\rightarrow\omega_{1}$ is upper semi-continuous, let $\lambda(r)=\sup\{\lambda(x)\mid x\in K\}$,
let $K_{\xi}(r)=\{x\in K\mid r(x)>\xi\}$ and let $\rho_{U}(r)=\min\{\xi\mid K_{\xi}(r)\text{ is compact}\}$.
As $r$ is upper semi-continuous, a routine compactness argument shows
that $\lambda(r)$, and hence $\rho_{U}(r)$, are countable ordinals.
If $K=\emptyset$, so that $r$ is the empty function, let $\lambda(r)=\rho_{U}(r)=0$
\begin{corollary}
\label{Cor uniform unique}
\begin{enumerate}
\item A uniform $\omega$-Stone space $W$ is determined up to isomorphism
by the triple $[K(W),r_{W},\rho(W)]$. 
\item If $[K,r,\rho]$ is such that $K$ is either $\mathscr{D}_{0}$ or
$\mathscr{D}_{1}$, $r\colon K\rightarrow\omega_{1}$ is upper semi-continuous
and $\rho\in\omega_{1}$ with $\rho_{U}(r)\leqslant\rho\leqslant\lambda(r)$,
then there is a uniform $\omega$-Stone space $W$ such that $(K(W),r_{W})\cong(K,r)$
and $\rho(W)=\rho$, and $W$ will be strongly uniform if $\rho_{U}(r)=\rho$.
\end{enumerate}
\end{corollary}
\begin{svmultproof}
If $W$ is a uniform $\omega$-Stone space, then by Theorem~\ref{Thm: uniform decomposn}
we can write $W=X\dotplus Y$, with $X$ strongly uniform, and $Y$
empty if $W$ is strongly uniform and otherwise scattered of type
$(\rho(W),\rho(W))$. Therefore $(K(X),r_{X})=(K(W),r_{W})$, and
the first statement follows from Theorem~\ref{Thm: strongly uniform unique}
and Theorem~\ref{Thm: Flum Martinez}. 

Conversely, given a tuple $[K,r,\rho]$ with the stated properties,
apply Theorem~\ref{Thm: strongly uniform unique} to find a strongly
uniform space $X$ such that $(K(X),r_{X})\cong(K,r)$. It follows
immediately that $\rho_{U}(X)=\rho_{U}(r)$ and $\text{\ensuremath{\nu(X)=}}\lambda(X)=\lambda(r)$.

If now $\rho=\rho_{U}(r)$, let $W=X$, so that $\rho(W)=\rho_{U}(W)=\rho$
as $W$ is strongly uniform. Otherwise let $W=X\dotplus Y$, where
$Y$ is scattered of type $(\rho,\rho)$; $W$ is uniform since $\nu(Y)=\rho\leqslant\lambda(r)=\nu(X)$.
\end{svmultproof}

\begin{theorem}
\label{Thm decomposn}An $\omega$-Stone space $W$ which is neither
uniform nor scattered can be written $W=X\dotplus Y$, where $X$
is strongly uniform of type $(K(W),r_{W})$ and $Y$ is:
\begin{enumerate}
\item \label{enu:type1}scattered of type $(\nu(W),0,n(W))$ if $\rho_{U}(W)=\rho(W)<\nu(W)$;
\item \label{enu:type2}scattered of type $(\nu(W),\nu(W))$ if $\rho_{U}(W)<\rho(W)=\nu(W)$;
\item \label{enu:type3}scattered of type $(\nu(W),\rho(W),n(W))$ if $\rho_{U}(W)<\rho(W)<\nu(W)$.
\end{enumerate}
\end{theorem}
\begin{svmultproof}
Let $\nu=\nu(W)$; the assumptions on $W$ mean that $\rho_{U}(W)<\nu$.
By decomposing $W$ as a disjoint union of compact opens, splitting
each such into its uniform and scattered components, recombining and
applying Corollary~\ref{Cor:dpcompact uniform}, we obtain $W=X\dotplus Y$,
where $X$ is strongly uniform with $\rho_{U}(X)=\rho_{U}(W)$, and
$Y$ is scattered. By Lemma~\ref{uniformaddition}, $\nu(X)<\nu(Y)$,
as $W$ is not uniform, and so $\nu(Y)=\nu$.

By Theorem~\ref{Thm: Flum Martinez} and Corollary~\ref{Ordspacetypes},
$Y$ has one of 3 types:
\begin{description}
\item [{Case~\ref{enu:type1}:}] $Y$ has type $(\nu,0,n)$, when $\rho_{U}(W)=\rho(W)$
as $Y$ is compact;
\item [{Case~\ref{enu:type2}:}] $Y$ has type $(\rho,\rho)$, when $\rho=\rho(W)=\nu$;
\item [{Case~\ref{enu:type3}:}] $Y$ has type $(\nu,\rho,n)$, where
$0<\rho<\nu$ and $\nu$ is a successor ordinal.
\end{description}
But in case~\ref{enu:type3}, we can write $Y=Y_{1}\dotplus Y_{2}$,
where $Y_{1}$ has type $(\nu,0,n)$ and $Y_{2}$ has type $(\rho,\rho)$,
and if $\rho\leqslant\rho_{U}(X)$ then $X\dotplus Y_{2}$ is also
strongly uniform, and we can replace $X$ with $X\dotplus Y_{2}$
and revert to case~\ref{enu:type1}. Therefore we only need to consider
case~\ref{enu:type3} for $\rho>\rho_{U}(W)$, when $\rho=\rho(W)$
as $\rho(X)=\rho_{U}(W)<\rho(Y)$. 

Finally, we have $n=n(W)$ in cases~\ref{enu:type1} and~\ref{enu:type3}.
\end{svmultproof}

Putting all these results together, we obtain:
\begin{theorem}[Uniqueness]

An $\omega$-Stone space $W$ is determined up to homeomorphism by
the tuple

\[
[(K(W),r_{W}),\nu(W),\rho(W),n(W)].
\]
\end{theorem}
\begin{svmultproof}
Immediate from Theorem~\ref{Thm: Flum Martinez} for scattered $W$,
Corollary~\ref{Cor uniform unique} for uniform $W$ and Theorem~\ref{Thm decomposn}
otherwise.
\end{svmultproof}

\begin{theorem}[Existence]

Suppose the tuple $[(K,r),\nu,\rho,n]$ is such that $K\in\{\emptyset,\mathscr{D}_{0},\mathscr{D}_{1}\}$,
$r\colon K\rightarrow\omega_{1}$ is upper semi-continuous, $n\in\mathbb{N}_{+}\cup\{-\infty\}$
and $\nu,\rho\in\omega_{1}$, with 
\begin{description}
\item [{(i)}] $\rho_{U}(r)\leqslant\rho$;
\item [{(ii)}] $\max(\rho,\lambda(r))\leqslant\nu$;
\item [{(iii)}] $\nu$ is a successor ordinal if $\rho<\nu$;
\item [{(iv)}] $n=-\infty$ precisely when $\max(\rho,\lambda(r))=\nu$.
\end{description}
Then there is an $\omega$-Stone space $W$ such that $[\nu(W),\rho(W),n(W)]=[\nu,\rho,n]$
and $(K(W),r_{W})\cong(K,r)$, and $W$ is unique up to homeomorphism.
\end{theorem}
\begin{svmultproof}
Suppose that a tuple $[(K,r),\nu,\rho,n]$ is given satisfying the
stated conditions. By Theorem~\ref{Thm: strongly uniform unique},
we can find a strongly uniform space $X$ such that $(K(X),r_{X})\cong(K,r)$,
with $\rho(X)=\rho_{U}(X)=\rho_{U}(r)$, letting $X=\emptyset$ if
$K=\emptyset$. 

If $\lambda(r)=\nu$, let $Y=\emptyset$ if $\rho=\rho_{U}(r)$ (which
includes the case when $\rho_{U}(r)=\rho=\nu$, as $\rho_{U}(r)\leqslant\lambda(r)$),
and let $Y$ be scattered of type $(\rho,\rho)$ if $\rho_{U}(r)<\rho\leqslant\nu$. 

Otherwise, if $\lambda(r)<\nu$, let $Y$ be scattered of type:

\[
\begin{cases}
(\nu,0,n) & \text{ if }\rho_{U}(r)=\rho<\nu\\
(\nu,\nu) & \text{ if }\rho_{U}(r)<\rho=\nu\\
(\nu,\rho,n) & \text{ if }\rho_{U}(r)<\rho<\nu
\end{cases}
\]

Finally let $W=X\dotplus Y$, which can be readily verified to have
the required invariants, with $n=-\infty$ precisely when $\lambda(r)=\nu$
or $\rho=\nu$.
\end{svmultproof}

\begin{corollary}[Decomposition]
\label{Cor:decomposn}If $W$ is an $\omega$-Stone space, then we
can write $W=X\dotplus Y$, with $X$ and $Y$ clopen in $W$, where:
\begin{enumerate}
\item either $X$ is empty, or $X$ is strongly uniform with $\lambda(X)=\lambda(W)$,
$K(X)=K(W)$ and $r_{X}=r_{W}$;
\item $Y$ is empty if $\rho(W)=\rho_{U}(W)$ and $\nu(W)=\lambda(W)$;
otherwise $Y$ is scattered with:
\begin{description}
\item [{(i)}] $n(Y)=n(W)$,
\item [{(ii)}] $\rho(Y)=\rho(W)$, unless $\rho_{U}(W)=\rho(W)$ when $\rho(Y)=0$;
\item [{(iii)}] $\nu(Y)=\nu(W)$, unless $\lambda(W)=\nu(W)$ when $\nu(Y)=\rho(W)$.
\end{description}
\end{enumerate}
\end{corollary}
\begin{svmultproof}
This is a straightforward consequence of Theorem~\ref{Thm: Flum Martinez}
for scattered $W$, Theorem~\ref{Thm: uniform decomposn} for uniform
$W$ and Theorem~\ref{Thm decomposn} otherwise, noting that if $Y$
is non-empty then $n(W)=n(Y)=-\infty$ precisely when $\nu(W)=\lambda(W)$
or $\nu(W)=\rho(W)$.
\end{svmultproof}

\section{\label{sec:Primitive-spaces,-trim}Primitive spaces, trim partitions
and their invariants}

So far we have considered general $\omega$-Stone spaces and explored
the invariants which determine their structure. \emph{Primitive spaces}
form an important sub-class of $\omega$-Stone spaces. For each primitive
space there is an associated canonical \emph{PO system} (poset with
a distinguished subset) which, together with an associated \emph{trim
partition} of the space, significantly determines the space's structure,
so that we can study primitive spaces by means of their associated
PO systems.

In this section we will extend the definition of some of the invariants
of Section~\ref{sec:Boolean-algebra-invariants} to cover PO systems,
and assess the extent to which the invariants for a primitive space
can be determined from those of an associated PO system. 

We refer to Pierce~\cite[section~3]{PierceMonk} for a detailed study
of the class of primitive Boolean algebras.

\subsection{Primitive Boolean rings and spaces}
\begin{definition}
If $R$ is a Boolean ring with Stone space $W$, we say that $A\in R$
is \emph{pseudo-indecomposable (PI)} if for all $B\in R$ such that
$B\subseteq A$, either $(B)\cong(A)$ or $(A-B)\cong(A)$; and that
$R$ and $W$ are \emph{primitive }if every element of $R$ is the
disjoint union of finitely many PI elements. We will say that $W$
is pseudo-indecomposable if it is compact and $1_{R}$ is PI\@.
\end{definition}
\begin{remark}
Hanf's original definition of primitive Boolean algebras in~\cite{Hanf}
required the unit of the algebra itself to be PI\@. We are following
most subsequent authors in dropping this requirement; we also no longer
require the Boolean ring to have a $1$. We note in passing that a
countable Boolean ring is primitive iff it is the direct product of
countable primitive Boolean algebras.
\end{remark}

\subsection{Trim partitions}
\begin{definition}
(PO systems) We recall that a \emph{PO system }is a set $P$ with
an anti-symmetric transitive relation $<$; equivalently, it is a
poset with a distinguished subset $P_{1}=\{p\in P\mid p<p\}$. We
write $p\leqslant q$ to mean $p<q$ or $p=q$; and $P^{d}$ for $\{p\in P\mid p\nless p\}$.

Let $P$ be a poset or PO system and let $Q\subseteq P$. We write
$Q_{\downarrow}=\{p\in P\mid p\leqslant q\text{ for some }q\in Q\}$,
$Q_{\uparrow}=\{p\in P\mid p\geqslant q\text{ for some }q\in Q\}$,
and write $P_{\min}$ and $P_{\max}$ for the sets of minimal and
maximal elements of $(P,\leqslant)$ respectively.

We recall that $Q$ is an \emph{upper} (respectively \emph{lower})
subset of $P$ if for all $q\in Q$, if $r\geqslant q$ (respectively
$r\leqslant q$) then $r\in Q$.

$P$ has the \emph{ascending chain condition (ACC)} if it has no strictly
increasing chain of elements $p_{1}\lneqq p_{2}\lneqq p_{3}\lneqq\cdots$.

A subset $Q$ of $P$ has a \emph{finite foundation} if there is a
finite subset $F\subseteq Q_{\downarrow}$ such that $Q_{\downarrow}\subseteq F_{\uparrow}$,
and has a \emph{finite ceiling }if there is a finite $F\subseteq Q$
such that $Q\subseteq F_{\downarrow}$.

An \emph{extended PO system }is a triple $(P,L,f)$, where $P$ is
a PO system, $L$ is a lower subset of $P_{\Delta}$, and $f\colon L_{\min}^{d}\rightarrow\mathbb{N}_{+}$,
where $P_{\Delta}=\{p\in P\mid\{p\}\text{ has a finite foundation}\}$
and $L_{\min}^{d}=L_{\min}\cap P^{d}$. 
\end{definition}
\textbf{Notation}: If $\{X_{p}\mid p\in P\}$ is a partition of a
subset $X$ of $W$, so that each $X_{p}\neq\emptyset$, we will write
$\{X_{p}\mid p\in P\}^{*}$ for the partition $\{\{X_{p}\mid p\in P\},W-X\}$
of $W$, with the convention that $W-X$ may be the empty set. The
asterisk will be omitted where a complete partition of $W$ is intended
(i.e.\ where $X=W$). 

For $Q\subseteq P$, we write $X_{Q}$ for $\bigcup_{p\in Q}X_{p}$.

We recall the following definitions from~\cite{Apps-Stone}.
\begin{definition}
(Partitions) Let $W$ be a Boolean space, $(P,<)$ a PO system, and
$\mathscr{X}=\{X_{p}\mid p\in P\}^{*}$ a partition of $W$; let $X=X_{P}$. 

We will say that $\mathscr{X}$ is a \emph{$P$-partition }of $W$
if $X_{p}^{\prime}\cap X=\bigcup_{q<p}X_{q}$ for all $p\in P$.

For a subset $Y$ of $W$, we define its \emph{type }$T(Y)=\{p\in P\mid Y\cap X_{p}\neq\emptyset\}$. 

A compact open subset $A$ of $W$ is \emph{$p$-trim,} and we write
$t(A)=p$, if $T(A)=\{q\in P\mid q\geqslant p\}$, with additionally
$|A\cap X_{p}|=1$ if $p\in P^{d}$. A set is \emph{trim }if it is
$p$-trim for some $p\in P$. We write $\Trim(\mathscr{X})$ for the
set of trim subsets of $W$ and $\Trim_{p}(\mathscr{X})$ for $\{A\in\Trim(\mathscr{X})\mid t(A)=p\}$. 

A $P$-partition of $W$ is \emph{complete} if $X=W$.

A $P$-partition $\mathscr{X}$ of $W$ is a \emph{trim }$P$-partition
if it also satisfies:
\begin{description}
\item [{\textbf{T1}}] Every element of $W$ has a neighbourhood base of
trim sets (this implies that $X$ is dense in $W$);
\item [{\textbf{T2}}] The partition is \emph{full}: that is, for each $p\in P$,
every element of $W$ with a neighbourhood base of $p$-trim sets
is an element of $X_{p}$;
\item [{\textbf{T3}}] For each $p\in P$, every element of $X_{p}$ has
a $p$-trim neighbourhood. 
\end{description}
We will say that $\mathscr{X}$ is a \emph{trim $(P,L,f)$-partition
of $W$ }if it is a trim $P$-partition of $W$ such that $\overline{X_{p}}$
is compact iff $p\in L$ and $|X_{p}|=f(p)$ for $p\in L_{\min}^{d}$,
where $L$ is a lower subset of $P$ and $f\colon L_{\min}^{d}\rightarrow\mathbb{N}_{+}$.
We note that $X_{p}$ is finite iff $p\in L_{\min}^{d}$, as by T1
and T2 $X_{p}$ is closed for $p\in P_{\min}$. Moreover, all (respectively
no) elements of $X_{p}$ are isolated in $X_{p}$ for $p\in P^{d}$
(respectively $p\notin P^{d}$).
\end{definition}
\begin{remark}
The structure of $P$ tells us little about the compactness structure
of $\mathscr{X}$ and nothing about the actual size of any finite
elements of $\mathscr{X}$. The above $(P,L,f)$ definition remedies
this by identifying which elements of the partition have a compact
closure, and the size of each finite such element. We showed in~\cite[Theorem~6.1]{Apps-Stone}
that the triple $(P,L,f)$ uniquely determines an $\omega$-Stone
space $W$ together with a $P$-partition of~$W$, subject to certain
``good behaviour'' conditions on $L$ and on $W$ itself.
\end{remark}
Trim partitions and primitive Boolean spaces are inextricably linked:
\begin{theorem}[{Apps~\cite[Theorem~3.12]{Apps-Stone}}]
\label{Primitive=00003Dtrim}An $\omega$-Stone space is primitive
iff it admits a trim $P$-partition for some PO system $P$.
\end{theorem}
In fact, a primitive $\omega$-Stone space $W$ has a canonical associated
PO system and trim partition. Namely, let $\mathscr{S}(W)$ be the
PO system comprising the homeomorphism classes of the compact open
pseudo-indecomposable subsets of $W$, with the relation $[A]<[B]$
iff $[A]\dotplus[B]\cong[A]$, where $[A]$ denotes the homeomorphism
class of a compact open subset $A$ of $W$. (This is the same as
the \emph{structure diagram }of the associated Boolean ring defined
by Hanf~\cite{Hanf}, but with the order relation reversed.) The
PO system $\mathscr{S}(W)$ is \emph{simple}: that is, any morphism
from it to another PO system is injective. The canonical trim partition
$\mathscr{X}(W)$ is then given by $\mathscr{X}(W)=\{X_{p}\mid p\in\mathscr{S}(W)\}^{*}$,
where $X_{p}$ is the set of points with a neighbourhood base of sets
homeomorphic to $p$, for $p\in\mathscr{S}(W)$. Indeed, the map $W\mapsto\mathscr{S}(W)$
induces a bijection between the homeomorphism classes of compact primitive
PI $\omega$-Stone spaces and the isomorphism classes of countable
simple PO systems with a unique minimal element (\cite[Theorem~5.5]{Hanf},~\cite[Corollary~3.8.3]{PierceMonk}).

In what follows, we will consider a general trim partition of a primitive
$\omega$-Stone space; any such partition is a refinement of the canonical
partition of the space, and the underlying PO system is the image
of the canonical PO system under some morphism (\cite[Corollary~3.21]{Apps-prim}).

\subsection{Trim partitions of scattered spaces}

We illustrate these ideas by describing the canonical trim partition
for a scattered space; for some purposes, this description of scattered
spaces may be easier to work with than the ordinals. 

\textbf{Notation}: If $\nu$ is an ordinal, let $(N(\nu),\prec)$
be the PO system consisting of the ordinals $\{\xi\mid\xi<\nu\}$,
with order reversal and $N(\nu)^{d}=N(\nu)$, so that $p\nprec p$
for all $p\in N(\nu)$ and $p\preccurlyeq q$ iff $p\geqslant q$
as ordinals. We note that $(N(\nu),\prec)$ has the ACC\@.

If $W$ is a scattered $\omega$-Stone space, let $P_{W}=N(\nu(W))$,
$L_{W}=\{\xi\in P_{W}\mid\xi\preccurlyeq\rho(W)\}$ and let $f_{W}\colon(L_{W})_{\min}^{d}\rightarrow\mathbb{N}_{+}$
be the empty function if $\rho(W)=\nu(W)$ (when $L_{W}=\emptyset$),
and otherwise let $f_{W}(\mu)=n(W)$ if $\nu(W)=\mu+1$, noting that
if $\rho(W)<\nu(W)$ then $\nu(W)$ is a successor ordinal and $(L_{W})_{\min}^{d}=\{\mu\}$.
\begin{theorem}
\label{Thm:scattered trim partition}Let $W$ be a scattered $\omega$-Stone
space. Then $W$ is primitive, and $\{W^{[\xi]}\mid\xi<\nu(W)\}$
is a complete trim $(P_{W},L_{W},f_{W})$-partition of~$W$.
\end{theorem}
\begin{svmultproof}
Let $\nu=\nu(W)$. By Lemma~\ref{Lemma: Wx-Wx+1}, $(W^{[\xi]})^{\prime}=\bigcup_{\eta>\xi}W^{[\eta]}$
and $\overline{W^{[\xi]}}=W^{(\xi)}$. Hence $\{W^{[\xi]}\mid\xi<\nu\}$
is a complete $P_{W}$-partition of $W$. To show that the partition
is trim, it suffices to show that if $w\in W^{[\xi]}$, $A\subseteq W$
is compact and open, and $w\in A$, then we can find a $\xi$-trim
set $B$ such that $w\in B\subseteq A$. To see this, we observe that
$W^{(\xi+1)}$ is a closed subset of $W$, so we can find a compact
open $B$ such that $w\in B\subseteq A\cap(W-W^{(\xi+1)})$ and $B\cap W^{[\xi]}=\{w\}$,
as $W^{[\xi]}$ is a discrete subset of $W-W^{(\xi+1)}$. So $B$
is $\xi$-trim as $w\in(W^{[\eta]})^{\prime}$ and $B\cap W^{[\eta]}\neq\emptyset$
for each $\eta\succ\xi$. 

Hence by Theorem~\ref{Primitive=00003Dtrim} $W$ is primitive.

Moreover, $\overline{W^{[\xi]}}$ is compact iff $\xi\geqslant\rho(W)$.
So if $\rho(W)=\nu(W)$ (which is always the case if $\nu(W)$ is
a limit ordinal), then none of the $\overline{W^{[\xi]}}$ are compact,
and $W$ admits a $(P_{W},L_{W},f_{W})$-partition where $L_{W}=\emptyset$
and $f_{W}$ is the empty function. Otherwise we have $\nu(W)=\mu+1$
with $W^{(\mu)}$ compact, and $n(W)=|W^{(\mu)}|$; thus $\{W^{[\xi]}\mid\xi\in P_{W}\}$
is a complete trim $(P_{W},L_{W},f_{W})$-partition of $W$.
\end{svmultproof}

\begin{remark}
$\{W^{[\xi]}\mid\xi<\nu(W)\}$ is the canonical trim partition of
$W$, as $W^{[0]}$ contains all the isolated points of $W$ (which
are their own neighbourhood base), and $W^{[\xi]}$ for $\xi>0$ contains
all points of $W$ with a neighbourhood base of sets homeomorphic
to $\omega^{\xi}+1$. 

The extended PO system $(P_{W},L_{W},f_{W})$ corresponds directly
to the invariants $(\nu(W),\rho(W),n(W))$ and so to the type $\tau(W)$
of $W$.

The primitivity of a scattered $\omega$-Stone space can also be shown
directly. Because an infinite compact open subspace of such a space
is homeomorphic to $n$ copies of $T(\omega^{\mu}+1)$ for some $\mu\in\omega_{1}$
and $n\in\mathbb{N}_{+}$, it is enough to show that $T(\omega^{\mu}+1)$
is pseudo-indecomposable. But by Theorem~\ref{Thm: Flum Martinez}
and~(\ref{ordCB}), any compact open subset of $T(\omega^{\mu}+1)$
containing $\omega^{\mu}$ is homeomorphic to $T(\omega^{\mu}+1)$. 
\end{remark}

\subsection{Trim partition preliminaries}

For this subsection, let $P$ be a PO system and let $\mathscr{X}=\{X_{p}\mid p\in P\}^{*}$
be a trim $P$-partition of the primitive $\omega$-Stone space $W$.

We will need the next two Lemmas, which show how the restriction of
$\mathscr{X}$ to certain subsets of $W$ remains trim. For $w\in W$,
let $V_{w}=\{A\in\Trim(\mathscr{X})\mid w\in A\}$ and let $I_{w}=\{t(A)\mid A\in V_{w}\}$,
being the trim neighbourhoods of~$w$ and their types.
\begin{lemma}
\label{Lemma:trim subsets}
\begin{enumerate}
\item \label{enu:trim1}(cf~\cite[Lemma~1.11(2)]{Hansoul2}) If $A$ is
a clopen subset of $W$, then $\{X_{p}\cap A\mid p\in T(A)\}^{*}$
is a trim $T(A)$-partition of $A$;
\item \label{enu:trim2}If $L$ is a lower subset of $P$, then $\{X_{p}\mid p\in L\}^{*}$
is a trim $L$-partition of $\overline{X_{L}}$.
\end{enumerate}
\end{lemma}
\begin{svmultproof}
\ref{enu:trim1}~is an easy check, using the fact that $(X_{p}\cap A)^{\prime}=X_{p}^{\prime}\cap A$.

\ref{enu:trim2}~Let $C=\overline{X_{L}}$ and $\mathscr{Y}=\{X_{p}\mid p\in L\}^{*}$,
being a partition of $C$. Let $R$ be the ring of compact opens of
$W$. Since $X_{L}\cap X_{q}^{\prime}=\bigcup_{r<q}X_{r}$ for $q\in L$,
we see that $\mathscr{Y}$ is an $L$-partition of $C$. If $x\in C$
and $A$ is a trim neighbourhood of $x$ in $W$, then $A\cap X_{L}\neq\emptyset$,
so $t(A)\in L$ and it follows easily that $A\cap C$ is a $t(A)$-trim
neighbourhood of $x$ in~$C$. Trim axioms T1 and T3 follow immediately
for $C$. Suppose finally that $x\in C$ has a neighbourhood base
in $C$ of $p$-trim sets, and let $D$ be one such. Then $V_{x}$
is a neighbourhood base for $x$ in $W$, and if $A\in V_{x}$ then
$A\cap C$ is a $t(A)$-trim neighbourhood of $x$ in $C$. But the
compact opens in $C$ have form $C\cap B$ for $B\in R$, and so $\{A\cap C\mid A\in V_{x}\}$
is a neighbourhood base of $x$ in $C$. It follows that $t(A)=p$
for all $A\in V_{x}$ such that $A\cap C\subseteq D$; hence $x$
has a neighbourhood base of $p$-trim sets in $W$, and $x\in X_{p}$
as required.
\end{svmultproof}

\begin{lemma}
\label{closures}Let $Q$ be a subset of $P$. Then $\overline{X_{Q}}=\{w\in W\mid I_{w}\subseteq Q_{\downarrow}\}$.
In particular, $\overline{X_{Q}}\cap X_{P}=X_{Q_{\downarrow}}$, and
$X_{Q}$ is closed iff $Q$ is a lower subset of $P$ satisfying the
ACC\@.
\end{lemma}
\begin{svmultproof}
For $w\in W$ such that $I_{w}\subseteq Q_{\downarrow}$ and $A\in V_{w}$,
we have $t(A)\in I_{w}\subseteq Q_{\downarrow}$ and $A\cap X_{Q}\neq0$.
Hence $\overline{X_{Q}}\supseteq\{w\in W\mid I_{w}\subseteq Q_{\downarrow}\}$.
Conversely, suppose that $w\in\overline{X_{Q}}$. For each $A\in V_{w}$,
$A\cap X_{Q}\neq0$, so we can find $q\in Q$ such that $q\geqslant t(A)$;
hence $t(A)\in Q_{\downarrow}$, and so $I_{w}\subseteq Q_{\downarrow}$.
Therefore $\overline{X_{Q}}=\{w\in W\mid I_{w}\subseteq Q_{\downarrow}\}$. 

If $w\in X_{p}$, then $p\in I_{w}$, so $w\in\overline{X_{Q}}$ iff
$p\in Q_{\downarrow}$; hence $\overline{X_{Q}}\cap X_{P}=X_{Q_{\downarrow}}$.
So $X_{Q}$ is closed only if $Q$ is a lower subset of $P$. By Lemma~\ref{Lemma:trim subsets}\ref{enu:trim2},
$\{X_{q}\mid q\in Q\}^{*}$ is a trim $Q$-partition of $\overline{X_{Q}}$
if $Q$ is a lower subset, and so $X_{Q}$ is closed iff $\{X_{q}\mid q\in Q\}$
is a complete trim $Q$-partition of $\overline{X_{Q}}$. But by~\cite[Theorem~5.1(a)]{Apps-Stone},
$\{X_{q}\mid q\in Q\}$ is a complete trim $Q$-partition of $\overline{X_{Q}}$
iff $Q$ has the ACC, as required.
\end{svmultproof}

We also recall the following result about isolated points:
\begin{proposition}[{\cite[Proposition~2.15]{Apps-Stone}}]
\label{Prop: isol}The isolated points of $W$ are precisely $X_{Q}$,
where $Q=P^{d}\cap P_{\max}$.
\end{proposition}

\subsection{PO system invariants}

We now mimic the Boolean space invariant definitions for a PO system,
and in Section~\ref{subsec:Primitive-space-invariants} show that
the Cantor-Bendixson sequence for a primitive $\omega$-Stone space
can be derived from a trim partition of the space together with the
Cantor-Bendixson sequence for the PO system underlying the trim partition.
The definitions below for $\lambda(P)$ and $r_{P}$ are elucidated
by Lemma~\ref{closures} above.
\begin{definition}
Let $P$ be a countable PO system. Then the \emph{Cantor-Bendixson
derivative }of $P$ is defined as $P^{\prime}=P-(P_{\max}\cap P^{d})$,
and the \emph{Cantor-Bendixson sequence} of $P$ is the list 

\[
(P^{(0)},P^{(1)},\ldots,P^{(\xi)},\ldots),
\]

where $P^{(0)}=P$, $P^{(\xi+1)}=(P^{(\xi)})^{\prime}$, and $P^{(\eta)}=\bigcap_{\xi<\eta}P^{(\xi)}$
if $\eta$ is a limit ordinal. Write $P^{[\xi]}=P^{(\xi)}-P^{(\xi+1)}=P_{\max}^{(\xi)}\cap P^{d}$.

Further, let

\begin{eqnarray*}
\nu(P) & = & \min\{\xi\mid P^{(\xi)}=P^{(\xi+1)}\};\\
K(P) & = & P^{(\nu(P))}\text{, the \emph{perfect} \emph{kernel} of }P;\\
\lambda(P) & = & \min\{\xi\mid P^{(\xi)}-K(P)\text{ is a lower subset of }P\};\\
r_{P}(p) & = & \min\{\xi\mid p\notin(P^{(\xi)}-K(P))_{\downarrow}\}\text{ for }p\in K(P);\\
K_{\xi}(P) & = & \{p\in K(P)\mid r_{P}(p)>\xi\}\text{, which is a lower subset of }K(P).
\end{eqnarray*}

The mapping $r_{P}\colon K(P)\rightarrow\omega_{1}$ is the \emph{rank
function} of $P$, and is an order reversing function.

$P$ is \emph{scattered }if $K(P)=\emptyset$ and is \emph{uniform}
if $\nu(P)=\lambda(P)$.
\end{definition}
\begin{remark}
As we shall see below, the analogues of $n(W)$ and $\rho(W)$ will
be provided by $f$ and $L$ respectively within an extended PO system
$(P,L,f)$, provided $P$ is ``reduced''.
\end{remark}
\begin{proposition}
\label{Prop: PO sys CB}Let $P$ be a countable PO system.
\begin{enumerate}
\item \label{enu:CB1}$P^{(\xi)}\supseteq P-P^{d}$ and $P-K(P)\subseteq P^{d}$;
\item \label{enu:CB2}$\nu(P)\geqslant\lambda(P)$;
\item \label{enu:CB3}each $P^{(\xi)}$ and $K(P)$ are lower subsets of
$P$, and $P-K(P)$ an upper subset of $P$;
\item \label{enu:CB4}$P^{(\xi)}-K(P)$ satisfies the ACC for each $\xi$;
\item \label{enu:CB5}$\lambda(P)=\sup\{r_{P}(p)\mid p\in K(P)\}$;
\item \label{enu:CB6}$(P^{[\xi]})_{\downarrow}=(P^{(\xi)}-K(P))_{\downarrow}=(P^{(\xi)}-K(P))\cup K_{\xi}(P)$.
\end{enumerate}
\end{proposition}
\begin{svmultproof}
\ref{enu:CB1}~is clear by induction on $\xi$, and~\ref{enu:CB2}
is immediate.

For~\ref{enu:CB3}, $P^{\prime}$ is a lower subset of $P$ as $P-P^{\prime}\subseteq P_{\max}$,
so by induction each $P^{(\xi)}$ is also a lower subset of $P$.

For~\ref{enu:CB4}, suppose $p_{1}\lneqq p_{2}\lneqq p_{3}\lneqq\cdots$
is a strictly increasing sequence in $P^{(\xi)}-K(P)$. As $P^{(\xi)}-K(P)=\bigcup\{P^{[\eta]}\mid\xi\leqslant\eta<\nu(P)\}$,
we can find $\{\eta_{n}\mid n\geqslant1\}$ such that $p_{n}\in P^{[\eta_{n}]}$
and $\xi\leqslant\eta_{n}<\nu(P)$ for each $n$. As $p_{n-1}\lneqq p_{n}\in(P^{(\eta_{n})})_{\max}$
and $P^{(\eta_{n})}$ is a lower subset of $P$, we have $\eta_{n-1}>\eta_{n}$.
So $\{\eta_{n}\}$ is a decreasing sequence of ordinals, which is
a contradiction. Hence $P^{(\xi)}-K(P)$ satisfies the ACC\@.

For~\ref{enu:CB5}, $P^{(\lambda(P))}-K(P)$ is a lower subset of
$P$ and hence $r_{P}(p)\leqslant\lambda(P)$ for all $p\in K(P)$.
However, if $\xi<\lambda(P)$, then $P^{(\xi)}-K(P)$ is not a lower
subset of $P$, so we can find $p\notin P^{(\xi)}-K(P)$ and $q\in P^{(\xi)}-K(P)$
such that $p\lneqq q$. Hence $p\in P^{(\xi)}$ as $P^{(\xi)}$ is
a lower subset, so we must have $p\in K(P)$ and $r_{P}(p)>\xi$,
and the result follows.

For~\ref{enu:CB6}, $r_{P}(p)>\xi$ iff $p\in(P^{(\xi)}-K(P))_{\downarrow}$
for $p\in K(P)$, and $P^{(\xi)}$ is a lower subset of $P$, so $(P^{(\xi)}-K(P))_{\downarrow}=(P^{(\xi)}-K(P))\cup K_{\xi}(P)$.
Clearly $(P^{[\xi]})_{\downarrow}\subseteq(P^{(\xi)}-K(P))_{\downarrow}$.
For equality, it is enough to show that $(P^{[\xi]})_{\downarrow}\supseteq P^{(\xi)}-K(P)$.
But if $p\in P^{(\xi)}-K(P)$, use~\ref{enu:CB4} to find a maximal
$q\in P^{(\xi)}-K(P)$ such that $p\leqslant q$. As $P^{(\xi)}-K(P)\subseteq P^{d}$
and $K(P)$ is a lower subset of $P$, we see that $q\in P_{\max}^{(\xi)}\cap P^{d}=P^{[\xi]}$,
as required.
\end{svmultproof}

\begin{remark}
It follows from~\ref{enu:CB6} that $P^{(\xi)}-K(P)$ is a lower
subset of $P$ iff $K_{\xi}(P)=\emptyset$. We therefore have the
following alternative formulae, as per Corollary~\ref{Cor: alt defns}: 

\begin{eqnarray*}
\nu(P) & = & \min\{\xi\mid P^{[\xi]}=\emptyset\};\\
\lambda(P) & = & \min\{\xi\mid K_{\xi}(P)=\emptyset\}=\min\{\xi\mid(P^{[\xi]})_{\downarrow}\cap K(P)=\emptyset\};\\
r_{P}(p) & = & \min\{\xi\mid p\notin(P^{[\xi]})_{\downarrow}\}\text{ for }p\in K(P)
\end{eqnarray*}
\end{remark}
Importantly, these invariants are preserved under surjective morphisms
(we recall that a \emph{morphism }$\alpha\colon Q\rightarrow P$ between
PO systems $Q$ and $P$ is a map $\alpha$ such that for each $q\in Q$,
$\{r\in Q\mid r>q\}\alpha=\{p\in P\mid p>q\alpha\}$):
\begin{proposition}
Let $P$ and $Q$ be PO systems and $\alpha\colon Q\rightarrow P$
a surjective morphism. Then $P^{[\xi]}=Q^{[\xi]}\alpha$ and $P^{(\xi)}=Q^{(\xi)}\alpha$
for all $\xi$, $\nu(P)=\nu(Q)$, $\lambda(P)=\lambda(Q)$, $K(P)=K(Q)\alpha$
and $r_{P}(q\alpha)=r_{Q}(q)$ for all $q\in K(Q)$.
\end{proposition}
\begin{svmultproof}
For $p\in P$, let $p^{+}=\{r\in P\mid r>p\}$, and similarly for
$q^{+}$ if $q\in Q$, so that $(q\alpha)^{+}=(q^{+})\alpha$. Now
$P^{[0]}=P^{d}\cap P_{\max}=\{p\in P\mid p^{+}=\emptyset\}$. Hence
$q\in Q^{[0]}$ iff $q\alpha\in P^{[0]}$, and $q\in Q^{\prime}=Q-Q^{[0]}$
iff $q\alpha\in P^{\prime}$. Therefore $Q^{\prime}\alpha=P^{\prime}$
and $\alpha$ restricts to a surjective morphism from $Q^{\prime}$
to $P^{\prime}$. It follows readily by transfinite induction that
for all $\xi$, $q\in Q^{[\xi]}$ iff $q\alpha\in P^{[\xi]}$, $q\in Q^{(\xi)}$
iff $q\alpha\in P^{(\xi)}$ and $Q^{(\xi)}\alpha=P^{(\xi)}$. For
example, if this is true for all $\xi<\eta$, then $Q^{(\eta)}\alpha=(\bigcap_{\xi<\eta}Q^{(\xi)})\alpha\subseteq P^{(\eta)}$,
while if $p\in P^{(\eta)}$ and $p=q\alpha$, say, then for all $\xi<\eta$
we have $q\alpha\in P^{(\xi)}$, so $q\in Q^{(\xi)}$ and $p\in Q^{(\eta)}\alpha$. 

It follows immediately that $\nu(P)=\nu(Q)$ and $K(P)=K(Q)\alpha$.
It remains to show that $r_{P}(q\alpha)=r_{Q}(q)$ for all $q\in K(Q)$,
for then the result for $\lambda(P)$ follows from Proposition~\ref{Prop: PO sys CB}\ref{enu:CB5}.
Now for $q\in K(Q)$ and $\xi\in\omega_{1}$, we have $q\alpha\notin P^{[\xi]}$,
and so $q\alpha\notin(P^{[\xi]})_{\downarrow}$ iff $(q\alpha)^{+}\cap P^{[\xi]}=\emptyset$
iff $(q^{+})\alpha\cap P^{[\xi]}=\emptyset$ iff $(q^{+})\cap Q^{[\xi]}=\emptyset$,
by the above. Hence $r_{P}(q\alpha)=\min\{\xi\mid q\alpha\notin(P^{[\xi]})_{\downarrow}\}=\min\{\xi\mid(q^{+})\cap Q^{[\xi]}=\emptyset\}=r_{Q}(q)$,
as required.
\end{svmultproof}

\subsection{\label{subsec:Primitive-space-invariants}Primitive space invariants
matched to those of an underlying extended PO system}

We now have the necessary machinery to link the relevant invariants
of a primitive space with those of any associated PO system. In this
section $W$ will be a primitive $\omega$-Stone space that admits
a trim $(P,L,f)$-partition for some extended PO system $(P,L,f)$. 

In the next Theorem where $W$ is a primitive $\omega$-Stone space
with a trim partition $\mathscr{X}=\{X_{p}\mid p\in P\}^{*}$, we
write $\mathscr{X}_{\lim}=W-X_{P}$, being the set of ``limit points''
of $\mathscr{X}$.
\begin{theorem}
\label{POsystemderivatives}Let $P$ be a PO system and $\mathscr{X}=\{X_{p}\mid p\in P\}^{*}$
a trim $P$-partition of the primitive $\omega$-Stone space $W$.
Then
\begin{enumerate}
\item \label{enu:deriv1}for each ordinal $\xi$, $W^{[\xi]}=X_{P^{[\xi]}}$,
$W^{(\xi)}=X_{P^{(\xi)}}\cup\mathscr{X}_{\lim}$, and $\{X_{p}\mid p\in P^{(\xi)}\}^{*}$
is a trim $P^{(\xi)}$-partition of the $\omega$-Stone space $W^{(\xi)}$;
\item \label{enu:deriv2}$\nu(W)=\nu(P)$, $\lambda(W)=\lambda(P)$, $K(W)=X_{K(P)}\cup\mathscr{X}_{\lim}=\overline{X_{K(P)}}$,
$K_{\xi}(W)=\overline{X_{K_{\xi}(P)}}$, $r_{W}(x)=r_{P}(p)$ for
$p\in K(P)$ and $x\in X_{p}$, and $r_{W}(x)=\min\{r_{P}(p)\mid p\in I_{x}\}$
for $x\in\mathscr{X}_{\lim}$.
\end{enumerate}
\end{theorem}
\begin{svmultproof}
\ref{enu:deriv1}~We claim first that $\mathscr{X}_{\lim}\subseteq\overline{X_{K(P)}}$.
Suppose $w\in\mathscr{X}_{\lim}$; by Lemma~\ref{closures}, it is
enough to show that $I_{w}\subseteq K(P)$, as $K(P)$ is a lower
subset of $P$. For otherwise let $p$ be a maximal element of $I_{w}\cap(P-K(P))$
(using Proposition~\ref{Prop: PO sys CB}\ref{enu:CB4}) and choose
a $p$-trim $A$ containing $w$. As the partition is trim, $w$ has
a neighbourhood base of trim subsets of $A$, and these must all be
$p$-trim by our choice of $p$, as $P-K(P)$ is an upper subset.
Hence $w\in X_{p}$ (as the partition is full), which gives the required
contradiction.

Now each $W^{(\xi)}$ is a closed subset of $W$ and so is itself
an $\omega$-Stone space. By Proposition~\ref{Prop: isol}, the isolated
points of $W$ are precisely the set $X_{Q}$ where $Q=P_{\max}\cap P^{d}$;
hence $W^{\prime}=X_{P^{\prime}}\cup\mathscr{X}_{\lim}$. If now $W^{(\xi)}=X_{P^{(\xi)}}\cup\mathscr{X}_{\lim}$,
it follows from Lemma~\ref{Lemma:trim subsets} that $\{X_{p}\mid p\in P^{(\xi)}\}^{*}$
is a trim $P^{(\xi)}$-partition of $W^{(\xi)}$, as $P^{(\xi)}$
is a lower subset of~$P$ (Proposition~\ref{Prop: PO sys CB}) and
$\overline{X_{P^{(\xi)}}}=W^{(\xi)}$ by the above. Applying Proposition~\ref{Prop: isol}
to $W^{(\xi)}$ we see that $W^{(\xi+1)}=X_{P^{(\xi+1)}}\cup\mathscr{X}_{\lim}$.
Moreover, for a limit ordinal $\eta$, $W^{(\eta)}=X_{P^{(\eta)}}\cup\mathscr{X}_{\lim}$
if this is true for all ordinals $\xi<\eta$. The second and third
assertions now follow by transfinite induction. The first assertion,
that $W^{[\xi]}=X_{P^{[\xi]}}$, also follows as $X_{P^{[\xi]}}\cap\mathscr{X}_{\lim}=\emptyset$. 

\ref{enu:deriv2}~It is immediate that $\nu(W)=\nu(P)$ and $K(W)=X_{K(P)}\cup\mathscr{X}_{\lim}$.
But $\mathscr{X}_{\lim}\subseteq\overline{X_{K(P)}}$ by the above
and $K(W)$ is closed, hence also $K(W)=\overline{X_{K(P)}}$. Now

\begin{align*}
\lambda(W) & =\min\{\eta\mid W^{(\eta)}-K(W)\text{ is closed}\}\\
 & =\min\{\eta\mid X_{P^{(\eta)}-K(P)}\text{ is closed}\},
\end{align*}
and it follows from Lemma~\ref{closures} that $\lambda(W)=\lambda(P)$,
as each $P^{(\eta)}-K(P)$ has the ACC (Proposition~\ref{Prop: PO sys CB}). 

Next, using Proposition~\ref{Prop: invariant basics}\ref{enu:prop3},
$K_{\xi}(W)=\overline{W^{[\xi]}}\cap K(W)=\overline{X_{P^{[\xi]}}}\cap K(W)$,
and $(P^{[\xi]})_{\downarrow}\cap K(P)=K_{\xi}(P)$. So $\overline{X_{K_{\xi}(P)}}\subseteq K_{\xi}(W)$,
using Lemma~\ref{closures}. But if $w\in K_{\xi}(W)$, then $I_{w}\subseteq K_{\xi}(P)$
by Lemma~\ref{closures}, and so $w\in\overline{X_{K_{\xi}(P)}}$
as every trim neighbourhood of $w$ meets $X_{K_{\xi}(P)}$.

If now $p\in K(P)$ and $x\in X_{p}$, then $r_{W}(x)=\min\{\xi\mid x\notin\overline{W^{[\xi]}}\}=\min\{\xi\mid p\notin(P^{[\xi]})_{\downarrow}\}=r_{P}(p)$,
by Lemma~\ref{closures}. Finally, for a general $x\in K(W)$, $x\notin\overline{W^{[\xi]}}$
iff $I_{x}\nsubseteq(P^{[\xi]})_{\downarrow}$ (Lemma~\ref{closures})
iff there is some $p\in I_{x}$ such that $p\notin(P^{[\xi]})_{\downarrow}$,
i.e. $r_{P}(p)\leqslant\xi$, and so $r_{W}(x)=\min\{r_{P}(p)\mid p\in I_{x}\}$
as required.
\end{svmultproof}

\begin{remark}
Let $W$ be a primitive $\omega$-Stone space admitting a trim $P$-partition
$\mathscr{X}$. By Corollary~\ref{Cor:decomposn}, we can write $W=Y\dotplus Z$
where $Y$ is strongly uniform or empty, and $Z$ is scattered or
empty. Restricting $\mathscr{X}$ to $Y$ and $Z$ in turn, we obtain
a trim $P_{u}$-partition $\mathscr{Y}$ of the strongly uniform space
$Y$ and a complete trim $P_{s}$-partition $\mathscr{Z}$ of the
scattered space $Z$, where $P=P_{u}\cup P_{s}$, $K(P)\subseteq P_{u}$
and $P_{s}\cap K(P)=\emptyset$. If also $P$ is neither uniform nor
scattered, and $|P^{[\xi]}|=1$ for each $\xi<\nu(P)$ (e.g.\ if
$P$ is reduced: see below), then $P_{s}\cap P_{u}=\bigcup\{P^{[\xi]}\mid\xi<\lambda(P)\}$.
So unlike Theorem~\ref{(Ketonen)}, $P$ does not generally split
into disjoint uniform and scattered subsets.
\end{remark}
We have seen that the invariants $\nu(W)$, $\lambda(W)$, $K(W)$
and the rank function $r_{W}$ can be derived directly from the corresponding
invariants for the PO system $P$ together with some limited information
about the partition (e.g.\ $r_{W}(x)$ depends only on $I_{x}$ and
$r_{P}$). We now consider $n(W)$ and $\rho(W)$, for which we will
need to make use of $f$ and $L$ respectively. Particularly for $\rho(W)$,
it will be convenient to restrict to ``reduced'' PO systems.

\textbf{Notation}: if $P$ is a PO system, we write $P_{0}$ for $\{p\in K(P)\mid r_{P}(p)=0\}$,
noting that $P_{0}=P-\overline{P^{[0]}}=K(P)-K_{0}(P)$.
\begin{definition}[Reduced PO systems]
A PO system $P$ is \emph{reduced }if $|P^{[\xi]}|=1$ for each $\xi<\nu(P)$
and $|P_{0}|\leqslant1$.
\end{definition}
\begin{remark}
If $P$ is a PO system, we can define the \emph{reduction of $P$}
as the PO system $\red(P)=K_{0}(P)\cup\widetilde{P_{0}}\cup N(\nu(P))$,
where $\widetilde{P_{0}}=\emptyset$ if $P_{0}=\emptyset$ and otherwise
$\widetilde{P_{0}}=\{\top\}$, with the following relations together
with those inherited from $(K_{0}(P),<)$ and $(N(\nu(P)),\prec)$:
\begin{enumerate}
\item $p<\xi$ if $p\in K_{\xi}(P)$ (recalling that the order relation
on $N(\nu(P))$ is reversed from that for normal ordinals);
\item $\top<\top$; 
\item for $p\in K_{0}(P)$, $p<\top$ iff $p<q$ for some $q\in P_{0}$.
\end{enumerate}
We can further define a natural surjective morphism $\red\colon P\rightarrow\red(P)$
by $\red(p)=p$ if $p\in K_{0}(P)$, $\red(p)=\top$ if $p\in P_{0}$,
and $\red(p)=\xi$ if $p\in P^{[\xi]}$. Hence the \emph{simple image}
of $P$, which is the image of $P$ after factoring out its largest
congruence, is reduced.

If now $\mathscr{X}=\{X_{p}\mid p\in P\}^{*}$ is a trim $P$-partition
of the $\omega$-Stone space $W$, then we can define the partition
$\red(\mathscr{X})=\{Y_{p}\mid p\in\red(P)\}^{*}$ of $W$, given
by $Y_{p}=X_{p}$ for $p\in K_{0}(P)$, $Y_{\xi}=W^{[\xi]}$ for $\xi\in N(\nu(P))$
and $Y_{\top}=K(W)-K_{0}(W)$ if $P_{0}\neq\emptyset$. It can be
shown that $\red(\mathscr{X})$ is a trim $\red(P)$-partition of
$W$, so we lose no information by assuming that $P$ is reduced.
\end{remark}
\begin{proposition}
\label{Prop:PLF formulae}Let $(P,L,f)$ be an extended PO system
and let $W$ be an $\omega$-Stone space that admits a trim $(P,L,f)$-partition. 
\begin{enumerate}
\item \label{enu:plf1}If $\nu(P)=\mu+1>\lambda(P)$ and $P^{[\mu]}$ is
a finite subset of $L$, then $n(W)=\sum\{f(p)\mid p\in P^{[\mu]}\}$
(and if $P$ is reduced, then $n(W)=f(p_{\mu})$, where $P^{[\mu]}=\{p_{\mu}\}$);
otherwise $n(W)=-\infty$;
\item \label{enu:plf2}if $P$ is a reduced PO system, then
\begin{description}
\item [{(i)}] $\rho(W)=\min\{\xi\mid P^{[\xi]}\subseteq L\}$;
\item [{(ii)}] $\min\{\xi\mid K_{\xi}(P)\subseteq L\}\leqslant\rho_{U}(W)\leqslant\min\{\xi\mid P^{[\xi]}\subseteq L\}$;
\item [{(iii)}] either $\rho_{U}(W)$ is a limit ordinal or $K_{\eta+1}(P)-K_{\eta}(P)$
is not a finite subset of~$L$, where $\rho_{U}(W)=\eta+1$.
\end{description}
\end{enumerate}
\end{proposition}
\begin{svmultproof}
\ref{enu:plf1}~If $\nu(P)=\mu+1>\lambda(P)$, then $W^{[\mu]}=\bigcup\{X_{p}\mid p\in P^{[\mu]}\}$,
using Theorem~\ref{POsystemderivatives}, which is finite iff $P^{[\mu]}$
is a finite subset of $L$, as $W^{[\mu]}$ is closed and discrete.
The result now follows from the definition of $n(W)$.

\ref{enu:plf2}~For (i), if $P$ is reduced, then $W^{[\xi]}=X_{P^{[\xi]}}$,
whose closure is compact iff $P^{[\xi]}\subseteq L$. For (ii), we
observe that $\rho_{U}(W)\leqslant\rho(W)$, and that if $K_{\xi}(P)\nsubseteq L$
then $K_{\xi}(W)$ is not compact as $K_{\xi}(W)=\overline{X_{K_{\xi}(P)}}$
(Proposition~\ref{POsystemderivatives}). For (iii), if $\rho_{U}(W)=\eta+1$
and $K_{\eta+1}(P)-K_{\eta}(P)$ were a finite subset of $L$ then
$\overline{K_{\eta}(W)}$ would be the union of $\overline{K_{\eta+1}(W)}$
and a finite set of compact sets and so would be compact.
\end{svmultproof}

\subsection{Uniqueness results for trim $(P,L,f)$-partitions}

We have seen that if $P$ is reduced and $W$ admits a trim $(P,L,f)$-partition,
then we can deduce $\nu(W)$, $\lambda(W)$, $\rho(W)$, and $n(W)$
from $(P,L,f)$, and can also largely determine $r_{W}$. However,
the question remains as to how the partition elements fit together,
and specifically to what extent the compact elements fit together
compactly. The issue comes into sharp focus when we consider $\rho_{U}(W)$,
for which we only have a range (per Proposition~\ref{Prop:PLF formulae})
rather than a formula.

We therefore ask the question, under what conditions on $(P,L,f)$
is there a unique $\omega$-Stone space admitting a trim $(P,L,f)$-partition?
We first recall some existing results:
\begin{theorem}
\label{Thm: uniqueness}Let $(P,L,f)$ be a countable extended PO
system.
\begin{enumerate}
\item \label{enu:uniq1}\cite[Theorem~6.1]{Apps-Stone}~If $L$ has a finite
foundation, then there is an $\omega$-Stone space~$W$ admitting
a bounded trim $(P,L,f)$-partition, unique up to $P$-homeomorphism.
\item \cite[Corollary~6.4]{Apps-Stone}~If $L$ has a finite ceiling, then
there is an $\omega$-Stone space $W$ with a trim $(P,L,f)$-partition,
unique up to $P$-homeomorphism, and the partition will be bounded.
\end{enumerate}
\end{theorem}
Here, a $(P,L,f)$-partition $\{X_{p}\mid p\in P\}^{*}$ is \emph{bounded
}if $\overline{X_{L}}$ is compact, and a homeomorphism between two
Boolean spaces admitting a $P$-partition is a \emph{$P$-homeomorphism
}if it respects the partitions.

It follows that if $L$ is finite or empty then there is a unique
$\omega$-Stone space $W$ with a trim $(P,L,f)$-partition. In particular,
for any countable $P$, there is a unique $\omega$-Stone space with
a trim $P$-partition $\{X_{p}\mid p\in P\}^{*}$ such that none of
the $\overline{X_{p}}$ are compact.

At the other end of the spectrum, we consider the case when $L$ is
almost as large as~$P$. Suppose that $P$ is a reduced PO system
and that $\nu(P)$ is a successor ordinal $\mu+1$. Suppose also that
$K_{0}(P)$ has a finite foundation $F$. This is equivalent to saying
that $P$ has a finite foundation, namely $F\cup P^{[\mu]}$, together
with $P_{0}$ in case $P_{0}$ is a minimal element. We observe also
that $P^{[0]}\cup P_{0}$ is a finite ceiling for $P$.

Firstly, let $L=P$ and suppose $f\colon P_{\min}^{d}\rightarrow\mathbb{N}_{+}$
is given. Then there is a unique (up to $P$-homeomorphism) $\omega$-Stone
space $W_{1}$ admitting a bounded $(P,P,f)$-partition, and this
space will be compact as $L=P$.

Secondly, if $P_{0}\neq\emptyset$, let $L=P-P_{0}$: this is a lower
subset of $P$, with a finite foundation $F\cup P^{[\mu]}$ and finite
ceiling $P^{[0]}$. So by the above we obtain a unique $\omega$-Stone
space $W_{2}$ admitting a $(P,P-P_{0},f)$-partition, and it is not
difficult to see that $W_{2}\cong W_{1}\dotplus\mathscr{D}_{0}$.

Thirdly, if $P^{[0]}\nsubseteq L$, then we may lose uniqueness if
$L$ does not have a finite ceiling: indeed, we may not even be able
to determine whether or not the space has a compact perfect kernel.
\begin{example}
\label{exa: non-unique PLf}Let $P=\{a_{n},b_{n},c\mid n\geqslant0\}$
with $P^{d}=\{a_{n},c\mid n\geqslant0\}$, and with the order relation
generated by the following:

$a_{n}>a_{n+1}$ and $a_{n}>b_{n}>c$ for all $n\geqslant0$.

We see easily that $\nu(P)=\omega$, with $P^{[n]}=\{a_{n}\}$ for
$n\geqslant0$ and $P_{0}=\emptyset$, so that $P$ is reduced, and
indeed simple. Hence $K(P)=\{b_{n},c\mid n\geqslant0\}$. Let $L=K(P)$,
which has a finite foundation $\{c\}$, and let $f(c)=1$. Applying~\cite[Theorem~5.5]{Apps-Stone}
with $Q=L$, taking $L^{u}=\emptyset$ and $L^{u}=L$ in turn, we
obtain $\omega$-Stone spaces $X$ and $Y$ with $(P,L,f)$-partitions
$\{X_{p}\mid p\in P\}^{*}$ and $\{Y_{p}\mid p\in P\}^{*}$ respectively
such that $\overline{X_{L}}$ is compact but $\overline{Y_{L}}$ is
not. But $K(X)=\overline{X_{L}}$ and $K(Y)=\overline{Y_{L}}$ by
Theorem~\ref{POsystemderivatives}. So $K(X)$ is compact and $K(Y)$
is not, and $X$ and $Y$ are not homeomorphic. With reference to
Theorem~\ref{Thm: uniqueness}\ref{enu:uniq1}, $X$ is bounded but
$Y$ is not; uniqueness applies only if we restrict to bounded spaces.

This example confirms that if $W$ admits a trim $(P,L,f)$-partition,
then we may not be able to deduce $\rho_{U}(W)$ from $(P,L,f)$.
For in this example $\rho_{U}(X)=0$ whereas $\rho_{U}(Y)=\omega$,
corresponding to the extremes of the possible range of values for
$\rho_{U}$ in Proposition~\ref{Prop:PLF formulae}.
\end{example}
\begin{remark}
This construction is a simpler version of that used by Pierce~\cite[Proposition~9.12]{Pierce}
to show that for any finite PO system $Q$, there is a finite simple
PO system $P$ such that $K_{0}(P)\cong Q$.
\end{remark}
If $P$ is not reduced but has a finite foundation, then an $\omega$-Stone
space $W$ with a trim $(P,P,f)$-partition need not even be compact.
For example, let $P=\{b_{n},c\mid n\geqslant1\}$, with $P^{d}=\{c\}$
and $b_{n}>c$ for all $n$, so that $P=K(P)$. Then $P$ has a finite
foundation, and applying~\cite[Theorem~5.5]{Apps-Stone} with $Q=L=P$
and $f(c)=1$, taking $L^{u}=\emptyset$ and $L^{u}=L$ in turn, we
obtain $\omega$-Stone spaces $X$ and $Y$ admitting $(P,P,f)$-partitions
such that $X$ is compact but $Y$ is not.

\section{\label{sec:Uniform-Boolean-rings}Uniform $\omega$-Stone spaces
and primitivity}

We have seen that a primitive $\omega$-Stone space $W$ admits a
trim $(P,L,f)$-partition for a suitable PO system $P$, and the invariants
derived from $P$ provide significant information about $W$.

In addition, all scattered $\omega$-Stone spaces are primitive. The
picture is more complex however for uniform $\omega$-Stone spaces.
We need only consider strongly uniform spaces (Theorem~\ref{Thm: uniform decomposn}),
and there is a 1--1 map between the homeomorphism classes of strongly
uniform spaces and isomorphism classes of pairs $(K,r)$, where $K\in\{\mathscr{D}_{0},\mathscr{D}_{1}\}$
and $r\colon K\rightarrow\omega_{1}$ is upper semi-continuous (Theorem~\ref{Thm: strongly uniform unique}).
An obvious question therefore is how to determine whether or not such
a space is primitive given its rank function. In this section we develop
two alternative conditions for a uniform $\omega$-Stone space to
be primitive based on the measure associated with the space's rank
function, which is another invariant of the space.

We recall the following definitions (e.g.\ from~\cite[section~1.12]{PierceMonk})
which link the rank function to a measure:
\begin{definition}
An \emph{m-monoid} is a commutative monoid whose zero element is the
unique unit.

Write $\mathscr{W}=\omega_{1}\cup\{o\}$, the set of all countable
ordinals together with a new zero element $o$ such that $o<\xi$
for all $\xi\in\omega_{1}$, which becomes an m-monoid with zero element
$o$ and binary operation $a+b=\max\{a,b\}$. 

If $R$ is a Boolean ring and $M$ is an m-monoid, an \emph{$M$-measure
on $R$ }is a map $\sigma\colon R\rightarrow M$ such that 
\begin{enumerate}
\item $\sigma(A\dotplus B)=\sigma(A)+\sigma(B)$ for all pairs $(A,B)$
of disjoint elements of $R$,
\item $\sigma(A)=o$ iff $A=0$.
\end{enumerate}
If $W$ is the perfect $\omega$-Stone space of the countable atomless
Boolean ring $R$, then upper semi-continuous maps $r\colon W\rightarrow\omega_{1}$
can be related to $\mathscr{W}$-measures $\sigma$ on $R$ as follows:
\begin{enumerate}
\item Given $r\colon W\rightarrow\omega_{1}$, define a measure $\sigma_{r}\colon R\rightarrow\mathscr{W}$
by $\sigma_{r}(0)=o$ and $\sigma_{r}(A)=\sup\{r(x)\mid x\in A\}$
if $A\in R-\{0\}$;
\item Given $\sigma\colon R\rightarrow\mathscr{W}$, define $r_{\sigma}(x)=\min\{\sigma(A)\mid x\in A\wedge A\in R\}$
for $x\in W$.
\end{enumerate}
\end{definition}
Let $R_{1}$ and $R_{0}$ denote the countable atomless Boolean rings
with and without a $1$ respectively. If $R=R_{1}$, then $\sigma_{r}$
is a $\mathscr{W}$-measure on $R$, $r_{\sigma}$ is an upper semi-continuous
map from the Stone space of~$R$ ($\mathscr{D}_{1}$) to $\omega_{1}$,
and the maps $r\mapsto\sigma_{r}$ and $\sigma\mapsto r_{\sigma}$
are inverse bijections (\cite[Proposition~1.12.2]{PierceMonk}). It
follows easily that this is also true if $R=R_{0}$, with Stone space
$\mathscr{D}_{0}$. The earlier uniqueness Theorem~\ref{Thm: strongly uniform unique}
translates into the following results for measures:
\begin{theorem}
\label{Thm: uniform<>rank function}\cite[Corollary~1.12.3]{PierceMonk}
For each $\mathscr{W}$-measure $\sigma\colon R_{1}\rightarrow\mathscr{W}$,
there is a compact uniform $\omega$-Stone space $X_{\sigma}$ such
that $(K(X_{\sigma}),r_{X_{\sigma}})\cong(\mathscr{D}_{1},r_{\sigma})$,
and the map $\sigma\mapsto X_{\sigma}$ gives a bijection between
isomorphism classes of $\mathscr{W}$-measures on $R_{1}$ and homeomorphism
classes of compact uniform $\omega$-Stone spaces.
\begin{theorem}
For each $\mathscr{W}$-measure $\sigma\colon R_{0}\rightarrow\mathscr{W}$,
there is a non-compact strongly uniform $\omega$-Stone space $X_{\sigma}$
such that $(K(X_{\sigma}),r_{X_{\sigma}})\cong(\mathscr{D}_{0},r_{\sigma})$,
and the map $\sigma\mapsto X_{\sigma}$ gives a bijection between
isomorphism classes of $\mathscr{W}$-measures on $R_{0}$ and homeomorphism
classes of non-compact strongly uniform $\omega$-Stone spaces.
\end{theorem}
\end{theorem}
We now have the building blocks for some further definitions.
\begin{definition}
Let $W$ be a Boolean space, and let $R=\Co(W)$ denote the Boolean
ring of compact open subsets of $W$. Let $M$ be an m-monoid and
$\sigma\colon R\rightarrow M$ an $M$-measure. For $A,B\in R$, write
$A\cong_{\sigma}B$ if there is an isomorphism $\alpha\colon(A)\rightarrow(B)$
such that $\sigma(C\alpha)=\sigma(C)$ for all $C\in(A)$. Then:
\begin{enumerate}
\item $A\in R$ is \emph{$\sigma$-self-similar at $x\in A\subseteq W$
}if there is a neighbourhood basis $S_{x}\subseteq R$ of $x$ such
that $B\cong_{\sigma}A$ for all $B\in S_{x}$;
\item $A$ is \emph{$\sigma$-pseudo-indecomposable ($\sigma$-PI)} if for
all $B\in R$ such that $B\subseteq A$, either $(B)\cong_{\sigma}(A)$
or $(A-B)\cong_{\sigma}(A)$;
\item $R$ and $W$ are \emph{$\sigma$-primitive} if every element of $R$
can be written as the disjoint union of finitely many $\sigma$-PI
elements of $R$.
\end{enumerate}
In the next Theorem, for a Boolean ring $R$ with Stone space $W$
we let $K(R)$ denote the Boolean ring of compact open subsets of
$K(W)$, so that $K(R)$ is trivial or countable atomless; elements
of $K(R)$ have the form $A\cap K(W)=K(A)$ for $A\in R$. Let $\sigma_{r}\colon K(R)\rightarrow\mathscr{W}$
denote the $\mathscr{W}$-measure corresponding to the rank function
$r\colon K(W)\rightarrow\omega_{1}$.
\end{definition}
\begin{theorem}
\label{Thm: sigma primitive}Let $W$ be an $\omega$-Stone space
and $R=\Co(W)$. Then $W$ is primitive iff $K(R)$ is $\sigma_{r}$-primitive.
\end{theorem}
\begin{svmultproof}
Suppose first that $W$ is primitive. We claim that if $A$ is PI
in $R$, then $K(A)$ is $\sigma_{r}$-PI in $K(R)$.

For suppose $K(A)=A\cap K(W)=F_{1}\dotplus F_{2}$, where $F_{i}=B_{i}\cap K(W)$
and $B_{i}\in R$ ($i=1,2$). We can replace $B_{1}$ with $A\cap B_{1}$
and $B_{2}$ with $A-B_{1}$ to assume that $A=B_{1}\dotplus B_{2}$.
But $A$ is PI, so (renumbering if necessary) there is an isomorphism
$\alpha\colon(B_{1})\rightarrow(A)$, so that $F_{1}\alpha=K(B_{1})\alpha=K(A)=A\cap K(W)$.
Moreover, $r(x\alpha)=r(x)$ for all $x\in F_{1}$, noting that $r_{B_{1}}(x)=r_{W}(x)$.
Hence $\sigma_{r}(G\alpha)=\sigma_{r}(G)$ for compact open subsets
$G$ of~$F_{1}$, so $F_{1}\cong_{\sigma_{r}}K(A)=A\cap K(W)$ and
$K(A)$ is $\sigma_{r}$-PI\@. 

If now $E\in K(R)$, with $E=K(A)$ where $A\in R$, write $A=A_{1}\dotplus\cdots\dotplus A_{n}$,
where each $A_{k}$ is PI\@. Then $E=K(A_{1})\dotplus\cdots\dotplus K(A_{n})$
with each $K(A_{k})$ being $\sigma_{r}$-PI\@. Hence $K(R)$ is
$\sigma_{r}$-primitive.

Now suppose that $K(R)$ is $\sigma_{r}$-primitive and $A\in R$.
We can write $A\cap K(W)=E_{1}\dotplus\cdots\dotplus E_{n}$, with
each $E_{k}$ being compact open and $\sigma_{r}$-PI in $K(R)$.
Find $B_{1},\ldots,B_{n}\in R$ such that $E_{k}=B_{k}\cap K(W)$
for $k\leqslant n$. Let $C_{1}=B_{1}\cap A$; for $2\leqslant k\leqslant n-1$,
let $C_{k}=B_{k}\cap A-\bigcup_{j\leqslant k-1}C_{j}$; and let $C_{n}=A-\bigcup_{j\leqslant n-1}C_{j}$;
then it is easy to see that $A=C_{1}\dotplus\cdots\dotplus C_{n}$
and that $C_{k}\cap K(W)=E_{k}$ for $k\leqslant n$.

Use Theorem~\ref{(Ketonen)} to write each $C_{k}=C_{k}^{u}\dotplus C_{k}^{s}$,
where $C_{k}^{u}$ is uniform or empty and $C_{k}^{s}$ is scattered
or empty. By Theorem~\ref{Thm:scattered trim partition} $C_{k}^{s}$
is primitive, and so can be written as a disjoint union of PI sets.
Clearly $C_{k}^{s}\cap K(W)=\emptyset$, and so $C_{k}^{u}\cap K(W)=E_{k}$.
We can now apply Lemma~\ref{Lemma: uniform pi} to obtain that $C_{k}^{u}$
is PI\@. Hence $A$ can be written as a disjoint union of PI sets,
as required.
\end{svmultproof}

\begin{lemma}
\label{Lemma: uniform pi}Let $W$ be an $\omega$-Stone space, $R=\Co(W)$,
and suppose $D\in R$ is such that $D$ is uniform and $D\cap K(W)$
is $\sigma_{r}$-PI, where $\sigma_{r}$ is the measure corresponding
to $r_{W}$. Then $D$ is pseudo-indecomposable.
\end{lemma}
\begin{svmultproof}
By Proposition~\ref{Prop: invariant basics}\ref{enu:prop5}, $\lambda(A)=\sigma_{r}(K(A))$
for $A\in R$. Suppose $D=D_{1}\dotplus D_{2}$, where each $D_{i}\in R$.
Then $K(D)=K(D_{1})\dotplus K(D_{2})$. As $K(D)$ is $\sigma_{r}$-PI,
we may suppose (renumbering if necessary) that $K(D)\cong_{\sigma_{r}}K(D_{1})$.
But $\nu(D_{1})\leqslant\nu(D)=\lambda(D)=\sigma_{r}(K(D))=\sigma_{r}(K(D_{1}))=\lambda(D_{1})\leqslant\nu(D_{1})$.
Hence $\lambda(D_{1})=\nu(D_{1})$, $D_{1}$ is uniform, and $D$
and $D_{1}$ have isomorphic rank functions as they are $\sigma_{r}$-isomorphic.
So by Theorem~\ref{Thm: invariants uniqueness} $D$ and $D_{1}$
are homeomorphic, and $D$ is pseudo-indecomposable.
\end{svmultproof}

We now give an alternative characterisation of $\sigma$-PI sets.
\begin{proposition}
\label{selfsimilar}Let $W$ be an $\omega$-Stone space, $R=\Co(W)$,
$M$ an m-monoid, $\sigma\colon R\rightarrow M$ an $M$-measure,
$A\in R$ and $x\in A\subseteq W$. The following are equivalent:
\begin{enumerate}
\item \label{enu:ss1}$A$ is \emph{$\sigma$-}self-similar at \emph{$x$,}
\item \label{enu:ss2}for all $B\in R$ such that $x\in B\subseteq A$,
we have $\sigma(B)=\sigma(A)$, and for all $A_{1}\subseteq A$ we
can find $B_{1}\subseteq B$ such that $A_{1}\cong_{\sigma}B_{1}$,
\item \label{enu:ss3}$B\cong_{\sigma}A$ for all $B\in R$ such that $x\in B\subseteq A$. 
\end{enumerate}
\end{proposition}
\begin{svmultproof}
\ref{enu:ss1}~$\Rightarrow$~\ref{enu:ss2}~Suppose $A$ is \emph{$\sigma$-}self-similar
at \emph{$x$ }with neighbourhood basis \emph{$S_{x}$. }Given $x\in B\subseteq A$
and $A_{1}\subseteq A$, find $E\in S_{x}$ such that $x\in E\subseteq B$
and $A\cong_{\sigma}E$. Then $\sigma(E)=\sigma(A)$, so by the additivity
of $\sigma$ we must have $\sigma(B)=\sigma(A)$. Let $\alpha\colon(A)\rightarrow(E)$
be a $\sigma$-preserving isomorphism and let $B_{1}=A_{1}\alpha$
to complete.

\ref{enu:ss2}~$\Rightarrow$~\ref{enu:ss3}~Suppose $x\in B\subseteq A$.
Define a relation between elements of $(A)$ and elements of $(B)$:
if $C\subseteq A$ and $D\subseteq B$, let $C\sim D$ if either $x\in C\cap D$,
or $x\notin C\cup D$ and $C\cong_{\sigma}D$. We must show that $\sim$
satisfies the Vaught conditions of Theorem~\ref{(Vaught)}. 

Clearly $A\sim B$, and $C\sim0$ iff $C=0$. Suppose that $C\sim D$,
and $D=D_{1}\dotplus D_{2}$. If $x\notin C\cup D$, then $C\cong_{\sigma}D$,
and we can immediately write $C=C_{1}\dotplus C_{2}$ with $C_{i}\cong_{\sigma}D_{i}$
for $i=1,2$. Suppose instead that $x\in C\cap D$ with $x\in D_{1}$.
Find $E\subseteq C$ such that $E\cong_{\sigma}D$, and write $E=E_{1}\dotplus E_{2}$
with $E_{i}\cong_{\sigma}D_{i}$ for $i=1,2$. 

Case 1: If $x\in C-E_{2}$, we can take $C_{1}=C-E_{2}$ and $C_{2}=E_{2}$. 

Case 2: If however $x\in E_{2}$, then we can find $F\subseteq D_{1}$
such that $E_{2}\cong_{\sigma}F$ and $G\subseteq E_{1}$ such that
$G\cong_{\sigma}F$, as $E_{1}\cong_{\sigma}D_{1}$. Hence $D_{2}\cong_{\sigma}E_{2}\cong_{\sigma}G$,
and $x\notin G$, so we can take $C_{1}=C-G$ and $C_{2}=G$ to complete.

The case when $C\sim D$, and $C=C_{1}\dotplus C_{2}$ is identical.

Hence by Theorem~\ref{(Vaught)} there is an isomorphism $\alpha\colon(A)\rightarrow(B)$
such that each $C\in(A)$ can be expressed as $C=C_{1}\dotplus\ldots\dotplus C_{n}$
where $C_{i}\sim C_{i}\alpha$ for all $i\leqslant n$. But then $\sigma(C_{i})=\sigma(C_{i}\alpha)$
as either $C_{i}\cong_{\sigma}C_{i}\alpha$, or $x\in C_{i}\cap C_{i}\alpha$
in which case $\sigma(C_{i})=\sigma(C_{i}\alpha)=\sigma(A)$. Hence
$\alpha$ is $\sigma$-preserving as $\sigma$ is additive, so that
$B\cong_{\sigma}A$ as required.

\ref{enu:ss3}~$\Rightarrow$~\ref{enu:ss1}~Immediate.
\end{svmultproof}

\begin{corollary}
\label{Cor:sigma pi}Let $W$ be an $\omega$-Stone space, $R=\Co(W)$,
$M$ an m-monoid, $\sigma\colon R\rightarrow M$ an $M$-measure,
and $A\in R$. 

Then $A$ is $\sigma$-pseudo-indecomposable iff $A$ is \emph{$\sigma$-}self-similar
at $x$ for some $x\in A$.
\end{corollary}
\begin{svmultproof}
$\Leftarrow$ Suppose $A=B\dotplus C$ (with $B,C\in R$), and say
$x\in B$. Then $B\cong_{\sigma}A$ by Proposition~\ref{selfsimilar}.
Hence $A$ is $\sigma$-PI\@.

$\Rightarrow$ Suppose that $A\in R$ is $\sigma$-PI and let $\{B_{1},B_{2},\ldots\}$
be an enumeration of $R$. Let $A_{0}=A$, and choose $A_{n}$ successively
for $n\geqslant1$ such that $A_{n}\subseteq A_{n-1}$, $A_{n}\cong_{\sigma}A$
and either $A_{n}\subseteq B_{n}$ or $A_{n}\cap B_{n}=\emptyset$,
by writing $A_{n-1}=(A_{n-1}\cap B_{n})\dotplus(A_{n-1}-B_{n})$ and
using inductively the fact that $A_{n-1}\cong_{\sigma}A$ and so $A_{n-1}$
is $\sigma$-PI\@. Then $\bigcap_{n\geqslant1}A_{n}$ is non-empty
by compactness and must be a singleton $\{x_{A}\}$, say, as for any
$x,y\in W$ we can find $n$ such that $x\in B_{n}$ and $y\notin B_{n}$.
Therefore $A$ is $\sigma$-self-similar at $x_{A}$. 
\end{svmultproof}

\begin{remark}
Taking $M=\mathscr{W}$, we obtain a convenient local topological
condition for testing when the rank function associated with a given
measure will yield PI sets in the associated uniform space, as we
will see in the final section.

Taking $M=\{o,0\}$ with the trivial $M$-measure, we recover the
result of Pierce~\cite[Proposition~3.1.1]{PierceMonk} that a compact
Boolean space $W$ is PI iff there is a point $x\in W$ such that
every clopen neighbourhood of $x$ is homeomorphic to $W$. Pierce
calls such an $x$ a \emph{point of homogeneity }for $W$; Pierce's
proof uses the fact that every PI Boolean space $W$ has the Schroder-Bernstein
property: namely that if $W$ and $X$ are each homeomorphic to a
compact open subset of the other, then they are homeomorphic to each
other.
\end{remark}

\section{\label{sec:A-non-primitive-countable}Non-primitive $\omega$-Stone
spaces: bad points}

To conclude this paper, we take a brief look at non-primitive Boolean
spaces.
\begin{definition}
Let $W$ be an $\omega$-Stone space. We will say that $w\in W$ is
a \emph{bad point }if it does not have a neighbourhood base of PI
sets.
\end{definition}
\begin{proposition}[{Pierce~\cite[Proposition~3.2.1]{PierceMonk}}]
An $\omega$-Stone space $W$ is non-primitive iff it contains a
bad point.
\end{proposition}
\begin{svmultproof}
We include an elementary proof of this result of Pierce. It is immediate
from the definition that a primitive space has no bad points. Suppose
instead that every element of $W$ has a neighbourhood base of PI
sets. By a standard compactness argument, it suffices to show that
a set which is a union of $n$ compact open PI subsets can be written
as a disjoint union of at most $n$ such sets. This is trivially true
for $n=1$. Suppose the statement is true for $n\leqslant k$, and
that $B=A_{1}\cup A_{2}\cup\cdots\cup A_{k}\cup A_{k+1}$, with each
$A_{i}$ compact open and PI in $W$. By the induction step, we may
assume that $\{A_{i}\mid i\leqslant k\}$ are disjoint or empty. 

We will adjust $A_{1}$ and/or $A_{k+1}$ so that we can apply the
induction hypothesis. As $A_{k+1}$ is PI, there are 2 cases, writing
$A_{k+1}=(A_{k+1}-A_{1})\dotplus(A_{k+1}\cap A_{1})$:

Case 1: if $A_{k+1}\cong A_{k+1}-A_{1}$, replace $A_{k+1}$ with
$A_{k+1}-A_{1}$ (which will equal $A_{k+1}$ if $A_{1}$ and $A_{k+1}$
are already disjoint). We now have $B=A_{1}\dotplus C$, where $C=A_{2}\cup\cdots\cup A_{k}\cup(A_{k+1}-A_{1})$,
and we can apply the induction hypothesis to $C$ to obtain $B$ as
the disjoint union of at most $k+1$ PI sets.

Case 2: if $A_{k+1}\cong A_{1}\cap A_{k+1}$, then $A_{1}=(A_{1}-A_{k+1})\dotplus(A_{1}\cap A_{k+1})\cong(A_{1}-A_{k+1})\dotplus A_{k+1}=A_{1}\cup A_{k+1}$.
We can therefore replace $A_{1}$ and $A_{k+1}$ with $A_{1}\cup A_{k+1}$,
which is PI, to obtain $B$ as the union of $k$ compact open PI sets,
and apply the induction hypothesis to obtain $B$ as the disjoint
union of at most $k$ such sets. 

Induction on $k$ now completes the proof.
\end{svmultproof}

\subsection{A construction method for Boolean spaces with bad points}

We describe a method for constructing $\mathscr{W}$-measures on $\mathscr{D}_{1}$
that will yield bad points in the associated Boolean spaces, which
will therefore be non-primitive. We illustrate this with a specific
example in Section~\ref{subsec:Construction}.
\begin{definition}
\label{def:completion}Let $\{\sigma_{n}\colon R_{n}\rightarrow\mathscr{W}\mid n\geqslant1\}$
be a set of $\mathscr{W}$-measures on countable atomless Boolean
algebras $R_{n}$, each with Stone space $X_{n}$ homeomorphic to
the Cantor set $\mathscr{D}_{1}$. Let $X$ be the one-point compactification
of $\{X_{n}\}$ at $c$, so that $X=\{c\}\cup\bigcup\{X_{n}\mid n\geqslant1\}$,
and let $R$ be the Boolean algebra of compact open subsets of $X$.
We view $R_{n}$ as a subring of $R$ and write $X_{n}$ for the multiplicative
identity of $R_{n}$. Let $A_{n}=\{c\}\cup\bigcup\{X_{m}\mid m\geqslant n\}$,
which is compact and open in~$X$, with $\{A_{n}\}$ forming a neighbourhood
basis for $c$. We define the \emph{completion }of $\{\sigma_{n}\}$
on $R$ to be the measure $\sigma\colon R\rightarrow\mathscr{W}$
defined as follows (noting that this gives a well-defined additive
measure):
\begin{enumerate}
\item for $B\in R_{n}$, $\sigma(B)=\sigma_{n}(B)$;
\item for $n\geqslant1$, $\sigma(A_{n})=\sup\{\sigma_{m}(X_{m})\mid m\geqslant n\}$.
\end{enumerate}
We will say that $\{(R_{n},\sigma_{n})\mid n\geqslant1\}$ are \emph{incompatible}
if for each $n$, there is no measure-preserving isomorphism between
$R_{n}$ and any principal ideal of $R_{-n}$, where $R_{-n}=\{A\in R\mid A\subseteq\bigcup_{m\neq n}X_{m}\}$,
and where a map $\alpha\colon R_{n}\rightarrow R$ is \emph{measure-preserving}
if $\sigma_{n}(A)=\sigma(A\alpha)$ for all $A\in R_{n}$.
\end{definition}
\begin{remark}
For the purposes of the next Theorem, the incompatibility definition
could be weakened to require that there is no measure-preserving isomorphism
between $R_{n}$ and any principal ideal of either (i) $R_{k}$, for
$k<n$ or (ii) $R_{>n}=\{A\in R\mid A\subseteq\bigcup_{m>n}X_{m}\}$,
\end{remark}
\begin{theorem}
With the notation above, the completion $(R,\sigma)$ of the incompatible
measures $\{(R_{n},\sigma_{n})\}$ is not $\sigma$-primitive.
\end{theorem}
\begin{svmultproof}
By Corollary~\ref{Cor:sigma pi} it is enough to show that if $D$
is a neighbourhood of $c$, then $D$ is not $\sigma$-self-similar
at any $x\in D$. Suppose instead that $D$ is self-similar at $x$. 

Case 1: $x\neq c$. Find $n$ such that $A_{n}\subseteq D-\{x\}$,
with $x\in X_{k}$, say, where $k<n$. Let $E\in R$ be any clopen
subset of $D\cap X_{k}$ such that $x\in E$. A $\sigma$-isomorphism
$\alpha\colon(D)\rightarrow(E)$ would restrict to give a measure-preserving
isomorphism between $R_{n}$ and the ideal $(X_{n}\alpha)$ of $R_{k}$,
which is not possible. Hence $D$ is not self-similar at $x$.

Case 2: $x=c$. Find $n$ such that $A_{n}\subseteq D$ and choose
any $E\in R$ such that $c\in E\subseteq A_{n+1}$. Suppose $\alpha\colon(D)\rightarrow(E)$
is a $\sigma$-isomorphism, and let $\beta\colon D\rightarrow E$
be the corresponding homeomorphism.

If $c\beta\neq c$, then $c\beta\in X_{k}$ for some $k>n$. As $\beta$
is continuous, we can find $m>k$ such that $X_{m}\beta\subseteq X_{k}$,
giving a measure-preserving isomorphism between $R_{m}$ and the ideal
$(X_{m}\alpha)$ of $R_{k}$, which is a contradiction.

If instead $c\beta=c$, then $X_{n}\beta\subseteq E-\{c\}$, so $\alpha$
restricts to a measure-preserving isomorphism between $R_{n}$ and
the principal ideal $(X_{n}\alpha)$ of $R_{-n}$, which is again
a contradiction.
\end{svmultproof}

\begin{corollary}
With the notation above, if the measures $\{(R_{n},\sigma_{n})\}$
are incompatible, then there is a compact uniform $\omega$-Stone
space~$W$ such that $(K(W),r_{W})\cong(X,r_{\sigma})$, and $W$
is not primitive.
\end{corollary}
\begin{svmultproof}
Immediate from Theorem~\ref{Thm: uniform<>rank function} and Theorem~\ref{Thm: sigma primitive}.
\end{svmultproof}

\begin{remark}
The constructed space $W$ will have a single bad point, namely the
``image'' of $c$ in $K(W)$.
\end{remark}

\subsection{\label{subsec:Construction}Construction of a set of incompatible
measures}

We construct a set of ``incompatible subsets'' of the Cantor set,
whose characteristic functions will yield an incompatible set of measures.

We recall that a \emph{closure algebra }is a Boolean algebra $S$
together with an additional closure operator $^{C}$ such that $x\subseteq x^{C}$,
$x^{CC}=x^{C}$, $(x\cup y)^{C}=x^{C}\cup y^{C}$ (all $x,y\in S$),
and $0^{C}=0$; we write $\overline{x}$ for $x^{C}$.

For $k\geqslant0$, let $(Q_{k},\prec)$ be the PO system with elements
$\{0,1,2,\ldots,k,k+2\}$ and relations $i\prec j$ iff $i\geqslant j+2$
or $i=j\leqslant k$, so that $Q_{k}^{d}=\{k+2\}$; we note that $(Q_{k},\preccurlyeq)$
corresponds to a finite quotient of the Rieger-Nishimura lattice.
Then $2^{Q_{k}}$ becomes a closure algebra with closure operator
$\overline{P}=P_{\downarrow}$ for $P\subseteq Q_{k}$. For $n\leqslant k$,
let $P_{n,k}$ be the subset of $Q_{k}$ consisting of $\{n,n+1,\ldots,k,k+2\}$.
\begin{proposition}
\label{Function hk}For each $n\geqslant1$ there is a first order
closure algebra function $h_{n}$ with the following properties, writing
$h_{n}(C,\mathscr{C})$ for $h_{n}(C)$ evaluated in $\mathscr{C}$
(where $C\in\mathscr{C})$:
\begin{enumerate}
\item \label{enu:hk1}if $A$ is a clopen element of $\mathscr{C}$, then
$h_{n}(C,\mathscr{C})\cap A=h_{n}(C\cap A,\mathscr{\mathscr{C}}[A])$
for all $C\in\mathscr{C}$ and $n\geqslant1$, where $\mathscr{C}[A]$
denotes the sub-closure algebra $\{B\in\mathscr{C}\mid B\subseteq A\}$;
\item \label{enu:hk2}if $\alpha\colon\mathscr{C}\rightarrow\mathscr{E}$
is an isomorphism of closure algebras, then $h_{n}(C,\mathscr{C})\alpha=h_{n}(C\alpha,\mathscr{E})$
for all $C\in\mathscr{C}$ and $n\geqslant1$;
\item \label{enu:hk3}if $\mathscr{C}=2^{Q_{k}}$ and $C=Q_{k}-\{0\}$,
then $h_{n}(C,\mathscr{C})=P_{n,k}\text{ for }n\leqslant k$, $h_{k+1}(C,\mathscr{C})=\{k+2\}$
and $h_{k+2}(C,\mathscr{C})=\emptyset$.
\end{enumerate}
\end{proposition}
\begin{svmultproof}
Let $h_{1}(C)=C$, and for $n>1$ inductively define $h_{n}$ as follows:

\begin{align*}
h_{2}(C) & =C\cap(\overline{1_{\mathscr{C}}-C})\\
h_{n+1}(C) & =h_{n}(C)\cap\overline{(h_{n-1}(C)-h_{n}(C))}
\end{align*}

Properties~\ref{enu:hk1} to~\ref{enu:hk3} are now routine checks,
as $\overline{D}\cap A=\overline{D\cap A}$ for $D\in\mathscr{C}$
if $A$ is clopen. 
\end{svmultproof}

\begin{remark}
It follows that if $A$ and $B$ are clopen elements of $\mathscr{C}$
such that $A\dotplus B=1_{\mathscr{C}}$, then $h_{n}(C,\mathscr{C})=h_{n}(C\cap A,\mathscr{\mathscr{C}}[A])\dotplus h_{n}(C\cap B,\mathscr{\mathscr{C}}[B])$
for all $C\in\mathscr{D}$ and $n\geqslant1$. 

If $C$ is closed, $h_{n}(C,\mathscr{C})$ is the new element introduced
at the $n$-th stage in determining the closure sub-algebra of $\mathscr{C}$
generated by $C$. In the notation of~\cite[Theorem~7.2]{Apps-Stone},
$h_{n}(C,\mathscr{C})=W-V_{n}$, where $W=1_{\mathscr{C}}$ and $A=W-C$.
\end{remark}
\begin{definition}
A closed subset $C$ of a topological space $Y$ is of \emph{type
$Q_{k}$ in $Y$} if the atoms of $\mathscr{C}$, the closure subalgebra
of $2^{Y}$ generated by $C$, form a complete trim $Q_{k}$-partition
of $Y$, with $C=Y-D$, where $D$ is the atom of $\mathscr{C}$ corresponding
to $0\in Q_{k}$, and with any discrete atoms of $\mathscr{C}$ being
singletons.
\end{definition}
It follows that $\mathscr{C}$ is isomorphic to $2^{Q_{k}}$. We recall
from~\cite[Theorem~7.3~ and~Lemma~7.1]{Apps-Stone} that the Cantor
set $\mathscr{D}_{1}$ has a closed subset of type \emph{$Q_{k}$
}for each $k\geqslant0$\emph{, }unique up to homeomorphisms of $\mathscr{D}_{1}$.
\begin{proposition}
\label{Prop: hn}If the closed subset $C$ is of type \emph{$Q_{k}$}
in the topological space $Y$, then:
\begin{enumerate}
\item $h_{n+1}(C,2^{Y})$ is non-empty and has no isolated points if $n<k$;
\item $h_{k+1}(C,2^{Y})$ is a singleton;
\item $h_{n+1}(C,2^{Y})=\emptyset$ if $n>k$.
\end{enumerate}
\end{proposition}
\begin{svmultproof}
We note that $h_{n+1}(C,2^{Y})=h_{n+1}(C,\mathscr{C})$, where $\mathscr{C}$
is the closure subalgebra of $2^{Y}$ generated by $C$. The statements
now follow easily from the isomorphism between $\mathscr{C}$ and
$2^{Q_{k}}$, using Proposition~\ref{Function hk} and applying Proposition~\ref{Prop: isol}
in relation to isolated points: for $P_{n+1,k}^{d}\cap(P_{n+1,k})_{\max}=\emptyset$
in $2^{Q_{k}}$ for $n<k$, and $(k+2)\in Q_{k}^{d}$, so that $h_{k+1}(C,\mathscr{C})$
is a discrete atom of $\mathscr{C}$ and therefore a singleton.
\end{svmultproof}

\begin{example}
\label{exa:incompat}We consider a set of pairs $\{(X_{k},C_{k})\mid k\geqslant1\}$,
where each $X_{k}\cong\mathscr{D}_{1}$ and $C_{k}$ is a closed subset
of $X_{k}$ of type $Q_{k}$. Let $R_{k}=\Co(X_{k})$ and $\sigma_{k}$
be the characteristic function on $C_{k}$: that is, for $B\in R_{k}$,
$\sigma_{k}(B)=1$ if $B\cap C_{k}\neq\emptyset$ in $X_{k}$ and
$\sigma_{k}(B)=0$ if $B\cap C_{k}=\emptyset$ and $B\neq0$, with
$\sigma_{k}(0)=o$. Let $R$ and $\sigma$ be constructed as in Definition~\ref{def:completion},
and view each $R_{k}$ as a subring of $R$. Here, the completion
of $\{\sigma_{k}\}$ is just the characteristic function of the closed
subset $C$ in the one-point compactification at $c$ of $\{X_{k}\}$,
where $C=\{c\}\cup\bigcup\{C_{k}\mid k\geqslant1\}$.

We claim that $\{(R_{k},\sigma_{k})\mid k\geqslant1\}$ is an incompatible
set of measures.

For suppose that there is a $\sigma$-preserving isomorphism $\alpha\colon R_{k}\rightarrow(B)$,
where $B\in R_{-k}$. Let $\beta\colon X_{k}\rightarrow B$ be the
corresponding homeomorphism, which yields an isomorphism between the
closure algebras $2^{X_{k}}$ and $2^{B}$. Let $J=\{m\in\mathbb{N}\mid B\cap X_{m}\neq\emptyset$\},
which is finite and excludes $k$, and let $C_{J}=\bigcup_{m\in J}C_{m}$
and $B_{m}=B\cap X_{m}$ for $m\in J$. By Lemma~\ref{Charfn-measure}
below, we have $C_{k}\beta=B\cap C_{J}$. We then obtain the following,
using Proposition~\ref{Function hk}\ref{enu:hk1} twice:

\begin{align*}
h_{k+1}(C_{k},2^{X_{k}})\beta & =h_{k+1}(C_{k}\beta,2^{B})\\
 & =h_{k+1}(C_{J}\cap B,2^{B})\\
 & =\dotplus\{h_{k+1}(C_{m}\cap B_{m},2^{B_{m}})\mid m\in J\}\\
 & =\dotplus\{h_{k+1}(C_{m},2^{X_{m}})\cap B_{m}\mid m\in J\}
\end{align*}

But $h_{k+1}(C_{k},2^{X_{k}})$ is a singleton, whereas $h_{k+1}(C_{m},2^{X_{m}})\cap B_{m}$
is empty or has no isolated points for $m\neq k$, using Proposition~\ref{Prop: hn}
and the fact that $B_{m}$ is clopen in $X_{m}$. Therefore no such
measure-preserving $\alpha$ exists, and the measures $\{(R_{n},\sigma_{n})\mid n\geqslant1\}$
are incompatible.
\end{example}
\begin{lemma}
\label{Charfn-measure}Let $W$ be a Boolean space, $R=\Co(W)$, and
let $\sigma$ be the characteristic function of the closed subset
$C$ of $W$. Suppose $D,E\in R$ and $\alpha\colon(D)\rightarrow(E)$
is an isomorphism, corresponding to a homeomorphism $\beta\colon D\rightarrow E$.
Then $\alpha$ is a $\sigma$-isomorphism iff $(D\cap C)\beta=E\cap C$.
\end{lemma}
\begin{svmultproof}
``If'' is clear. For ``only if'', it suffices to show that $(D-C)\beta=E-C$.
But if $A\in(D)$ and $A\cap C=\emptyset$, then $\sigma(A)=0$, so
$\sigma(A\alpha)=0$ and $A\alpha\cap C=\emptyset$. Hence $(D-C)\beta\subseteq E-C$,
and equality follows by considering $\alpha^{-1}$. 
\end{svmultproof}

\subsection{A PI set can contain bad points}

In the previous example there was a single bad point, and no compact
open set containing it was PI\@. The picture is not always that simple
however.
\begin{example}
Let $(R,\sigma)$ be the completion of the incompatible measures of
Example~\ref{exa:incompat}, and let $(S_{n},\tau_{n})$ be isomorphic
to $(R,\sigma)$ for each $n\geqslant1$, with $Y_{n}$ the Stone
space of $S_{n}$ and $c_{n}$ being its bad point. Let $(S,\tau)$
be the completion of $\{(S_{n},\tau_{n})\mid n\geqslant1\}$, whose
Stone space $Y$ is the one-point compactification at $d$, say, of
$\{Y_{n}\}$. It can be shown that if $A$ is a compact open neighbourhood
in $Y$ of $d$, then $A\cong_{\tau}1_{S}$, so that $1_{S}$ is PI
and $d$ is a clean point. (This follows from the fact that for $B\in R_{k}$,
either $(B,\sigma_{k}|_{B})\cong(R_{k},\sigma_{k})$ or $(1-B,\sigma_{k}|_{1-B})\cong(R_{k},\sigma_{k})$.)
So the PI set $1_{S}$ contains bad points $\{c_{n}\}$, whose limit
is the clean point $d$.
\end{example}

\section{Some areas for further study}

We suggest the following questions.

\paragraph{\textbf{Ketonen's Boolean hierarchy }}
\begin{enumerate}
\item Compact $\omega$-Stone spaces can be classified by means of derived
measures and the Boolean hierarchy (e.g.\ see~\cite[1.21]{PierceMonk}).
To what extent does this classification still work, and/or still yield
uniqueness, for locally compact $\omega$-Stone spaces?
\item For primitive $\omega$-Stone spaces admitting a trim $P$-partition
for a PO system $P$, is there an equivalent of derived measures and
the Boolean hierarchy for PO systems?
\end{enumerate}

\paragraph{\textbf{Trim partitions of primitive spaces}}
\begin{enumerate}
\item What restrictions on an extended PO system $(P,L,f)$ result in uniqueness
of any $\omega$-Stone space admitting a trim $(P,L,f)$-partition?
For example, would uniqueness apply if $P$ were simple and $\nu(P)$
was finite?
\end{enumerate}

\paragraph{\textbf{Non-primitive spaces}}
\begin{enumerate}
\item Let $W$ be the Cantor set $\mathscr{D}_{1}$ and let $R=\Co(W)$,
the countable atomless Boolean algebra. Let $I$ denote the set of
finite binary strings (including the empty string), and let $\{A_{i}\in R\mid i\in I\}$
be a binary tree of non-zero elements of $R$ that generate~$R$:
i.e.\ $A_{\emptyset}=R$, and for $i\in I$, we have $A_{i}=A_{i0}\dotplus A_{i1}$.
Define a measure $\sigma\colon R\rightarrow\{0,1,o\}$ as follows:
\begin{description}
\item [{(i)}] $\sigma(0)=o$;
\item [{(ii)}] $\sigma(R)=1$;
\item [{(iii)}] If $\sigma(A_{i})=0$ ($i\in I$), then $\sigma(A_{i0})=\sigma(A_{i1})=0$;
\item [{(iv)}] If $\sigma(A_{i})=1$ ($i\in I)$, then $\sigma(A_{i1})=1$
and $\sigma(A_{i0})$ is $0$ or $1$ with probability $0.5$.
\end{description}
\begin{conjecture}
For the measure $\sigma$ constructed above, $K(W)$ is $\sigma$-primitive
with probability $0$. (We surmise that this is a similar situation
to algebraic and transcendental numbers: i.e.\ that ``almost all''
countable Boolean algebras are non-primitive, whereas the majority
of easily constructed such algebras are primitive.)
\end{conjecture}
\item How extensive can the set of bad points be in a non-primitive space?
For example, can they be dense in the perfect kernel (noting that
all points not in the perfect kernel are clean points)?
\end{enumerate}

\end{document}